\documentclass[reqno, 11pt]{smfart}
\usepackage[bbgreekl]{mathbbol}
\usepackage{amsmath,amssymb,amsthm,amsfonts, amscd, url, enumerate, leftidx}
\usepackage{appendix}
\usepackage[numeric]{amsrefs}
\usepackage{hyperref}
\usepackage[mathscr]{eucal}
\usepackage{comment}
\usepackage{stmaryrd}
\usepackage{multirow}
\usepackage{tikz-cd}
\usepackage{fancyhdr}

\input xy
\xyoption{all}

\makeatletter
\def\@tocline#1#2#3#4#5#6#7{\relax
  \ifnum #1>\c@tocdepth % then omit
  \else
    \par \addpenalty\@secpenalty\addvspace{#2}%
    \begingroup \hyphenpenalty\@M
    \@ifempty{#4}{%
      \@tempdima\csname r@tocindent\number#1\endcsname\relax
    }{%
      \@tempdima#4\relax
    }%
    \parindent\z@ \leftskip#3\relax \advance\leftskip\@tempdima\relax
    \rightskip\@pnumwidth plus4em \parfillskip-\@pnumwidth
    #5\leavevmode\hskip-\@tempdima
      \ifcase #1
       \or\or \hskip 1em \or \hskip 2em \else \hskip 3em \fi%
      #6\nobreak\relax
    \dotfill\hbox to\@pnumwidth{\@tocpagenum{#7}}\par
    \nobreak
    \endgroup
  \fi}
\makeatother

\usepackage{xcolor}
\usepackage{graphics}
\usepackage{graphicx}

\usepackage{mathptmx} % rm & math
\usepackage[scaled=0.90]{helvet} % ss
\usepackage{courier} % tt
\normalfont
\usepackage[T1]{fontenc}
\newtheorem*{nthm}{Theorem}

\theoremstyle{remark}

\theoremstyle{plain}
\newtheorem{theorem}{Theorem}[subsection]
\newtheorem{lemma}[theorem]{Lemma}
\newtheorem{proposition}[theorem]{Proposition}
\newtheorem{proposition-definition}[theorem]{Proposition-Definition}

\theoremstyle{definition} 
\newtheorem{definition}[theorem]{Definition}

\theoremstyle{remark}

\newtheorem*{rmk}{Remark}

\numberwithin{equation}{section}

\newcommand{\zee}{\mathbb{Z}}
\newcommand{\Q}{\mathbb{Q}}
\newcommand{\C}{\mathbb{C}}
\newcommand{\R}{\mathbb{R}}
\renewcommand{\O}{\mathcal{O}}
\newcommand{\Spec}{\mathrm{Spec}}

\DeclareMathOperator{\ad}{ad}
\DeclareMathOperator{\spanVector}{span}
\DeclareMathOperator{\depth}{dp}
\DeclareMathOperator{\aut}{Aut}

\DeclareMathOperator{\covol}{covol \;}

\newcommand{\deghat}{\widehat{\text{deg}} \;}

\newcommand{\be}[1]{\begin{eqnarray} \label{#1}}

\newcommand{\ee}{\end{eqnarray}}

\newcommand{\tpoint}[1]{\subsubsection{#1}}

\newcommand{\spoint}{\subsubsection{}}

\numberwithin{equation}{section}

\DeclareRobustCommand\longtwoheadleftarrow{\twoheadleftarrow\joinrel\relbar}

\DeclareRobustCommand\longhookrightarrow{\lhook\joinrel\longrightarrow}

\begin{document}

\title{Infinite-rank Euclidean Lattices and Loop Groups}

\author{Mathieu Dutour}
\email{dutour@ualberta.ca}
\author{Manish M. Patnaik}
\email{patnaik@ualberta.ca}

\maketitle

\begin{center} \textit{ To Mathukumalli Venkata Subbarao on the centennial of his birth, in gratitude } \end{center}

\vspace{0.15in} 
\begin{abstract}
In this paper, we associate a family of infinite-rank pro-Euclidean lattices to elements of a formal loop group and a highest weight representation of the underlying affine Kac--Moody algebra. In the case that the element has a polynomial representative, we can prove our lattices are theta-finite in the sense of Bost, allowing us to attach to each of our lattices a well-defined theta-like function.    

\end{abstract}

\setcounter{tocdepth}{1}
\tableofcontents

\newcommand{\GL}{\mathrm{GL}}
\newcommand{\uL}{\underline{L}}
\newcommand{\rk}{\mathrm{rank}}
\renewcommand{\gg}{\widehat{G}}
\newcommand{\pp}{\widehat{P}}
\newcommand{\ind}{\mathrm{Ind}}
\newcommand{\pro}{\mathrm{Pro}}
\newcommand{\Hom}{\mathrm{Hom}}
\newcommand{\lat}{\mathbf{Lat}}
\newcommand{\rr}{\rightarrow}
\newcommand{\lr}{\leftarrow}
\newcommand{\mc}[1]{\mathcal{#1}}
\newcommand{\ov}[1]{\overline{#1}}
\newcommand{\wh}[1]{\widehat{#1}}
\newcommand{\ardeg}{\widehat{\mathrm{deg}}}
\newcommand{\hilb}{\mathrm{hilb}}
\newcommand{\ssum}{\mathbf{Sum}}
\newcommand{\bG}{\mathbf{G}}
\newcommand{\kom}{\mathbb{K}}
\newcommand{\ade}{\mathbb{A}}
\newcommand{\ide}{\mathbb{I}}
\newcommand{\bU}{\mathbf{U}}
\newcommand{\bH}{\mathbf{H}}
\newcommand{\fin}{\mathrm{fin}}
\renewcommand{\gg}{\widehat{G}}
\renewcommand{\P}{\mathscr{P}}
\newcommand{\hP}{\wh{\P}}
\newcommand{\Ad}{\mathrm{Ad}}
\newcommand{\St}{\mathrm{St}}
\newcommand{\gen}{\mathsf{g}}
\newcommand{\mf}[1]{\mathfrak{#1}}

\newcommand{\dd}{\mathbf{d}}

\newcommand{\rts}{\mathscr{R}}
\newcommand{\crts}{\mathscr{R}^{\vee}}
\newcommand{\alphav}{\alpha^{\vee}}
\newcommand{\la}{\langle}
\newcommand{\ra}{\rangle} 
\newcommand{\lb}{\llbracket}
\newcommand{\rb}{\rrbracket}

\newcommand{\wts}{\mathscr{P}}
\newcommand{\lev}{\mathrm{lev}}
\newcommand{\dep}{\mathrm{dep}}
\newcommand{\gaff}{\mf{g}^e(\wh{A})}
\newcommand{\av}{a^{\vee}}
\newcommand{\cc}{\mathbf{c}} 
\newcommand{\kil}{\mathrm{kil}}
\renewcommand{\v}{\mathbf{v}}
\newcommand{\haff}{\widehat{\mf{h}}}
 \newcommand{\comp}{\mathrm{comp}}

\newcommand{\rtsaff}{\widehat{\rts}}
\newcommand{\rtsaffim}{\widehat{\rts}_{\mathrm{im}}}
\newcommand{\rtsaffre}{\widehat{\rts}_{\mathrm{re}}}
\newcommand{\affpi}{\widehat{\Pi}}
\newcommand{\rtl}{\wh{\mathscr{Q}}}

\newcommand{\vectZ}{\overline{\mathrm{Vect}}_{\mathbb{Z}}}
\newcommand{\provectZ}{\mathrm{pro}\, \overline{\mathrm{Vect}}_{\mathbb{Z}}}
\newcommand{\indvectZ}{\mathrm{ind}\, \overline{\mathrm{Vect}}_{\mathbb{Z}}}

\newcommand{\bb}{\widehat{B}}
\newcommand{\ovX}{\ov{\mc{X}}}

\section{Introduction}

The aim of this paper is to study a connection between two worlds of infinite-dimensional objects. On the one hand, we have the theory of loop groups and the representations of infinite-dimensional (affine) Kac--Moody algebras, and, on the other, certain infinite-rank metrized lattices. Whereas the theory of loop groups (in the arithmetic form which we need here) was started by H. Garland in the late 1970s \cite{gar:loopGroups}, the theory of infinite-rank lattices is of more recent vintage, having been formally introduced only a few years ago by J.-B. Bost (\textit{cf.} \cite{bost:book}). An interesting feature of Bost's theory is the existence of (several) so-called theta invariants attached to such lattices. Among their applications, a study of these invariants (see \cite[\S 10.8]{bost:book}) underlies some recent advances in diophantine algebraization results. Our contribution here is to construct certain families of infinite-rank Euclidean lattices with  finite theta invariants starting from elements in a loop group, or rather certain arithmetic quotients of it by parabolic subgroups. We believe this may be the first step in a fruitful interaction between the theory of loop groups, automorphic forms on them, and the infinite theta invariants studied by Bost. In particular, as we shall try to explain in more detail in the next section, our desire to extend certain aspects of the Langlands--Shahidi program to loop groups over number fields was the starting point of this work.

\subsection{Function fields} 

 Although this work is concerned with the number field  $\Q$, we begin with the case of a  function field $F$ of a smooth projective curve $C$ over a finite field $\mathbb{F}_q$. Write $\ade_F$ for the ring of adeles of $F$ and $\ide_F$ for the group of ideles with $|\cdot|$ the idelic norm. Recall that $\ade_F$ is the restricted product $\prod'_{v \in |F|} F_v$ over the set of places $|F|$ of $F$, where almost all the components of this product lie in the integral subrings $\O_v \subset F_v.$ Writing $\GL_n(F)$ for the general linear group on an $n$-dimensional $F$-vector space, it is an old observation, often attributed to A. Weil, that the double coset space \be{weil-dc} \begin{array}{lll} \mc{X}_{\GL_n} & := & \kom \setminus \GL_n(\ade_F) / \GL_n(F), \end{array} \ee with $\kom:= \prod_v \, \GL_n(\O_v)$, parametrizes isomorphism classes of principal $\GL_n$- bundles on $C.$ For an adelic element $g \in \GL_n(\ade_F)$, let $\P_g$ be a representative in this isomorphism class and $\Ad(\P_g):= \P_g \times_{\GL_n} \mathfrak{gl}_n$ be the corresponding vector bundle attached to the adjoint representation of $\GL_n$ on its Lie algebra $\mathfrak{gl}_n$. In addition to its rank and degree, this vector bundle has two cohomological invariants, namely \be{h0:h1} \begin{array}{lllllll} h^0(\Ad(\P_g)) & := & \log_q \, \left \vert  \, H^0(C, \Ad(\P_g)) \right \vert  & \text{ and } & h^1(\Ad(\P_g)) & := & \log_q \, \left \vert  \, H^1(C, \Ad(\P_g)) \right \vert \end{array} . \ee 

\begin{rmk} A vector bundle on $C$ may be seen as a \emph{coherent system of lattices} (\textit{cf}. \cite[p.97]{weil:bnt}) $\underline{L}:= ( L_v )_{v \in |F|}$  where coherence signifies that for almost all $v$, these lattices are equal to a fixed `trivial' or `reference' lattice. Thus the numbers (\ref{h0:h1}) can equally be regarded as invariants of some coherent family of lattices $\uL_g$. \end{rmk}

\newcommand{\st}{\mathrm{St}}
\noindent We might ask how to compute these invariants directly from $g \in \GL_n(\ade_F)$.  As for $h^0$,  one has \be{h0:stab} \begin{array}{lllllll} h^0(\Ad(\P_g)) & = &  \log_q \left \vert \st_{\GL_n(F)}(g)  \right \vert, & \text{ where } & \st_{\GL_n(F)}(g) & := & g^{-1} \kom g \cap \GL_n(F) \end{array} \ee is the stabilizer of $g$ inside the symmetric space $\kom \setminus \GL_n(\ade_F)$ under the right action of $\GL_n(F)$. The size of this group plays an important role in constructing the natural (Tamagawa) measure on $\mc{X}_{\GL_n}$.

Returning to $h^1$, first recall that by Serre duality \cite[Chap III, \S 7]{hartshorne}, if $E$ is a vector bundle on $C$,  \be{}\label{} \begin{array}{lll} H^1(C, E) & \cong & H^0(C, E^{\vee} \otimes \omega_C) \end{array} \ee where $E^{\vee}$ is the dual vector bundle and $\omega_C$ is the canonical line bundle on the curve. The Riemann--Roch theorem (\textit{cf}. \cite[Ch. 4, Thm 1.3]{hartshorne} for line bundles or \cite[pp. 96-101]{weil:bnt} for the higher rank situation in the language of lattices) further asserts that we have \be{rr:curves} \begin{array}{lll} h^0(C, E) - h^0(C, E^{\vee} \otimes \omega_C) & = &  \deg(E) + n \, (1 - \gen) , \end{array} \ee where $\gen$ is the genus of the curve, $n$ is the rank of $E$, and $\deg(E) = \deg ( \wedge^n \,E )$ is the degree of $E$.  For example, we have \be{g:1} \begin{array}{llll} h^0(C, E) & = & h^0(C, E^{\vee})  + \deg(E) & \text{ when } \gen =1. \end{array} \ee

\noindent In the case $E=\Ad(\P_g)$, we  saw in \eqref{h0:stab} how to compute $h^0$ from $g$, and we thus obtain \be{h0:h1:g1} \begin{array}{llll} \log_q  \left \vert \st_{\GL_n(F)}(g) \right \vert & = & h^0(\Ad(\P_g)^{\vee}) +  \deg(\Ad(\P_g)) & \text{ when } \gen=1. \end{array} \ee This formula, expressing the size of stabilizers in the function field analogue of a locally symmetric space in essentially geometric terms, played an important role in our thinking. Note that there is also a version of this formula for genera $\gen \neq 1$.

An important variant used in this paper of the construction above replaced the right action of $\GL_n(F)$ in~\eqref{weil-dc} with the action by the group of upper triangular matrices $B_n(F).$ The Iwasawa decomposition yields \be{B:G} \begin{array}{lllll} \mc{X}_{B_n} & := & \kom \setminus \GL_n(\ade_F) / \GL_n(F) & \cong & \kom \cap B_n(\ade_F)  \setminus B_n(\ade_F) / B_n(F). \end{array} \ee Using the first description above, we may consider the natural map $\mc{X}_{B_n } \longrightarrow  \mc{X}_{\GL_n}$, the fibers of which are called \emph{reductions} of a given principal $\GL_n$-bundle to a $B_n$-bundle. If we pick an element $x \in \mc{X}_{B_n}$, the second description in (\ref{B:G}) allows us to construct a principal $B_n$-bundle to which we can associate a vector bundle $\Ad_{\mf{b}_n}(\P_x)$ by taking the adjoint action of $B_n$ on the Lie algebra $\mathfrak{b}_n$. This bundle is related to stabilizers in $\mc{X}_{B_n(F)}$ and a formula for these groups can again be computed-- the answer now involves $2 \rho$, the sum of the positive roots of the Lie algebra $\mf{gl}_n$.

Let us now move to the case of loop groups while remaining within the realm of function fields. Let $G$ now be some simple, finite dimensional algebraic group, \textit{e.g.} $G= SL_n$. From $G$, we construct an infinite rank group over a base field $k$ as follows. First consider $LG(k):= G\left(k((t))\right)$, the group with points in the field of formal Laurent series over $k$. There are two important modifications one needs to make to $LG(k)$ to obtain the actual object of interest to us: first one considers a non-trivial (both in the mathematical and colloquial sense) central extension of this group by $k^*$; and second, one forms a semi-direct product with the automorphism of `loop rotations'  $g(t) \mapsto g(\tau t)$ with $g(t) \in G\left(k((t))\right)$ and $\tau \in k^*.$  In practice one fixes some value of $\tau$ and considers the corresponding object, written as $\gg^{\tau}(k)$. It is not a group (just a subset of one), though for the purposes of this introduction, we often regard it as such. We refer to (\ref{gEtaTau}) for the precise definition, and remark here that the affine analogues of 'discrete' and `maximal compact' subgroups still act on this set. We also note here that there is an adelic analogue of $\gg^{\tau}(k)$ in which $k$ is replaced by $\mathbb{A}_F$ and $\tau$ by an idele of $F$. Once again the analogue of the maximal compact subgroup, to be denoted by $\wh{\kom},$ and discrete subgroup, to be denoted as $wh{G}(F)$, subgroups act on this set $\wh{G}^{\tau}(\ad_F)$.

There is a important distinction we need to draw involving loop rotations: from the point of view of arithmetic quotients, what is important is the relation between $\tau$ and the direction in which the completion of the loop group is taken, \textit{i.e.} in either the \emph{positive}~$t$ or \emph{negative} $t^{-1}$ direction. For example, Garland \cite[Thm. 19.3]{gar:loopGroups}, who works with groups completed in \emph{positive powers of $t$}, requires $\left \vert \tau \right \vert<1$. In fact \be{weil:dc-lg} \begin{array}{lll} X^{\tau}_{\wh{G}} & = & \wh{\kom} \setminus \gg^{\tau}(\ade_F)  / \gg(F) \end{array} \ee is essentially compact (modulo the central extension) but this same reduction theory does not apply when we have $|\tau|>1$ and the corresponding space behaves `very' infinite-dimensionally. Nonetheless, applications from automorphic forms on loop groups suggests that one needs to confront this complexity.

 Geometrically, there are two ways to think of (\ref{weil:dc-lg})-- either one can regard its elements as parametrizing certain infinite-rank bundles on $C$ with symmetry `group' $\gg^{\tau}$ or, following M. Kapranov \cite{kap:sdual} (see also \cite{pat:thesis} for a more group-theoretic approach), as parametrizing \emph{finite-rank} $G$-bundles (with additional data\footnote{In fact, formulating the `other data' is a bit complicated as it involves the (relative) \emph{second} Chern classes and their relations to the central extensions of loop groups-- it is actually the main difficulty in interpreting this space.})  on a certain ruled (affine) \emph{surface} $S_{\tau} \rr C$ attached to $C$ (the curve corresponding to $F$) and $\tau$. We mostly adopt the former point of view, but let us mention a few motivating features from the surface picture. First, the condition $|\tau|<1$ was given a beautiful interpretation by Kapranov: it states that normal bundle to the natural embedding $C \hookrightarrow S_{\tau}$ given by the zero section has \emph{negative} degree.  Second, the bundle obtained from $x \in \gg^{\tau}(\ade_F)$ does not naturally live on $S_{\tau}$ but rather on the punctured surface $S_{\tau}^o:= S_{\tau} \setminus C$. It can however be extended to a bundle on $S_{\tau}$ and the Iwasawa decomposition (and the choices required to make one) give us one preferred way to do this. Denote this assignment of $x$ to a bundle on $S_{\tau}$ as $x \mapsto \hP_x.$ In terms of the infinite-rank picture, we should think of $\hP_x$ as corresponding to a \emph{reduction} of the $\gg^{\tau}$-bundle to a $\pp^{\tau}$-bundle on the curve $C$, where $\pp^{\tau} \subset \gg^{\tau}$ is some loop analogue of a parabolic subgroup.  Such ($\tau$-twisted) $\pp$-bundles on the curve naturally arise from elements in \be{weil:dc-pg} \begin{array}{lll} \mc{X}^{\tau}_{\pp} & = & \wh{\kom} \setminus \gg^{\tau}(\ade_F)  / \pp(F). \end{array} \ee Note that $\pp^{\tau}$ is the semi-direct product of a (finite-dimensional) torus and a \emph{pro}-unipotent group, so the elements from (\ref{weil:dc-pg}) can be naturally thought to parametrize some `pro'-bundle whereas elements from the original set \eqref{weil:dc-lg} would correspond to an `ind-pro'-bundle. There is also a \emph{non-compact} variant: replace $\pp$ with the corresponding negative parabolic subgroup $\pp^-$ (it is no longer a pro-group, but naturally only an ind-group) while keeping $|\tau|<1$, \textit{i.e.}  we consider \be{weil:dc-pg:m} \begin{array}{lll} \mc{X}^{\tau}_{\pp^-} & = & \wh{\kom} \setminus \gg^{\tau}(\ade_F)  / \pp^-(F). \end{array} \ee

\noindent Suppose as above we consider the adjoint representations $\mf{\widehat{p}}^{\pm}$ of $\pp^{\pm}$ and then form the corresponding vector bundles (of infinite-rank)  $\Ad_{\mf{\widehat{p}}^{\pm}}(\hP_x)$. As a set, we can again identify the the cohomology of the infinite-rank bundle with the stabilizers in $\mc{X}^{\tau}_{\pp^{\pm}}$, \textit{i.e.} we have $ H^0(C, \Ad_{\mf{\widehat{p}}^{\pm}}(\hP_x)) = x^{-1} \wh{\kom} x \cap \pp^{\pm}(F), $ where the cohomology group is defined from a certain Cech resolution. 

In the case $\left \vert \tau \right \vert < 1$, H. Garland constructed (in the number field setting, but the same argument works for function fields) a natural measure on $\mc{X}^{\tau}_{\pp}$ (and actually, $\mc{X}^{\tau}_{\gg}$ when $|\tau| \ll 1$, \textit{cf}. \cite[Appendix A]{gar:MS4}). However, the groups which one might expect to arise in the function field analogue of Garland's construction, namely the groups $x^{-1} \wh{\kom} x \cap \pp(F)$, are actually infinite ! Resolving this paradox, some years ago, A. Braverman and D. Kazhdan \cite{bk:ecm}, observed that Garland's measure could be naturally interpreted in the function field setting as a \emph{regularization} of the cardinality of $\vert x^{-1} \wh{\kom} x \cap \pp(F) \vert,$  effected by taking the left hand side of (\ref{h0:h1:g1}) and replacing it by the right hand side of Riemann--Roch applied to $\Ad_{\wh{\mf{p}}}(\hP_x)$. Ignoring the infinite constant involving the genus and the rank of this bundle, the analogue of the right hand side of \eqref{h0:h1:g1} is given in terms of  $H^0(C, \Ad_{\wh{\mf{p}}}(\hP_x)^{\vee})$ and $\deg(\Ad_{\wh{\mf{p}}}(\hP_x))$. Now the former space is finite-dimensional and the latter can be regularized using a well-known normal-ordering process in affine Kac--Moody theory (\textit{i.e.} the procedure of defining $2\rho$ for loop algebras as the `sum of all positive roots' in an affine Lie algebra). All of this crucially relies on the fact that we have $\left \vert \tau \right \vert<1$ since moving to the dual bundle gives a context where we have $\left \vert \tau \right \vert>1$, \textit{i.e.} an \emph{amplesness} condition, ensuring enough vanishing that the group $H^0(C, \Ad_{\wh{\mf{p}}}(\hP_x)^{\vee})$ is finite. On the other hand, if we work with $\pp^-$ instead (still keeping the condition $\left \vert \tau \right \vert<1$), then in fact the group $x^{-1} \kom x \cap \pp^-(F)$ is finite, \textit{i.e.} $H^0(C, \Ad_{\wh{\mf{p}}^-}(\hP_x))$ is finite.  This is roughly\footnote{We are glossing over an issue of completion, as $\pp^-$ is an ind-group but our main theta-finitness is for pro-objects.} an analogue of the main theta-finiteness result of the present work. 

Finally, we mention that a version of this type of Riemann--Roch regularization turns out to be crucial to the recent study of Eisenstein series on loop groups over function fields \cite{gmp}. There, we are in a context in which we need to cut down a certain non-convergent sum (the `naive' Eisenstein series attached to $\pp^-$) by an infinite set of the same flavour as $H^0(C, \Ad_{\wh{\mf{p}}}(\hP_x))$ as well as compensate for this by adding in a factor related to the finite quantity $H^0(C, \Ad_{\wh{\mf{p}}^-}(\hP_x)).$ These observations underlie the extension of the Langlands--Shahidi method to loop groups over function fields and it was in trying to generalize this method to number fields that we were led to the ideas described in this work.

\subsection{Number fields}

Let us now turn to the setting of the number field $F= \Q$. By a well-known analogy, a vector bundle over $C$ is now replaced by either of the following equivalent notions: a Hermitian vector bundle on the arithmetic curve $\Spec(\zee)$ or a Euclidean lattice, \textit{i.e.} a free $\zee$-module $E_{\zee}$ in a finite-dimensional vector space $E$ equipped with a Euclidean norm $\left \Vert \cdot \right \Vert$. Such objects are parametrized by the double coset space (\textit{cf}. \cite[p.216]{gillet}) \be{} \begin{array}{lllll} \ovX_{GL_n} & = & \kom \setminus GL_n(\ade_Q) / GL_n(\Q) & \cong & O(n) \setminus GL_n(\R) / GL_n(\zee) \end{array} \ee where $O(n)$ is the orthogonal group and $\kom \subset \GL_n(\ade_{\Q})$ is the direct product $\prod_{p} \, GL_n(\zee_p) \times O(n)$. The second equality above is a consequence of the strong approximation theorem for $\Q$.  Given $g \in GL_n(\R)$ we denote the corresponding Euclidean lattice as $\ov{E}_g:= (E_{g},\left \Vert \cdot \right \Vert_g)$ where $E_g$ is a free $\zee$-module and $\left \Vert \cdot \right \Vert_g$ a Euclidean norm on $E_{g, \R}:= E_g \otimes_{\zee} \R$. To tighten the analogy with what we described earlier, one can also develop a theory of metrized principal bundles for a general group and form associated Hermitian bundles (see \cite{chambertLoir-tschinkel:article:arithmeticTorsors}), but we do not pursue this point here.

As remarked by Bost \cite[\S 3.1]{bost:book}, there are (at least) two natural analogues of $h^0(C, E)$.  The first, \be{}\label{} \begin{array}{lll} h^0_{\mathrm{Ar}}(\overline{E}) & := & \log \, \left \vert \{ v \in L \; \mid \; \left \Vert v \right \Vert \leqslant 1 \} \right \vert, \end{array} \ee is more traditional in the context of Arakelov geometry, whereas the second, namely the \emph{theta invariant} \be{theta:inv} \begin{array}{lll} h^0_{\vartheta}(\overline{E}) & = & \log \, \sideset{}{}\sum\limits_{v \in E} e^{ - \pi \left \Vert v \right \Vert^2}, \end{array} \ee is the one which we shall be concerned with here. An immediate consequence of the Poisson summation formula (over $\Q$) is that this latter invariant is a perfect analogue to \eqref{h0:h1:g1}, as we have \be{}  \begin{array}{lll} h^0_{\vartheta}(\ov{E}) - h^0_{\vartheta}(\ov{E}^{\vee}) & = & \wh{\mathrm{deg}}\, {E} , \end{array} \ee whereas the corresponding result for $ h^0_{\mathrm{Ar}}$ leads only to an asymptotic statement. Another nice feature of the theta invariants which do not hold for their Arakelov counterparts is a certain subaddivity phenomenon in exact sequences. Nonetheless, there are now precise comparisons in finite dimensions between these two invariants (see \cite[Chap. 3]{bost:book} and references therein for more details). As far as we know, it is only the theta invariants which have been extended to infinite-rank lattices and, for this reason, we focus on them. 

Formalizing the notion of an infinite-rank lattice can be done in several different ways leading to the notion of ind-Euclidean lattices and pro-Euclidean lattices. \footnote{there should also be a notion of ind-pro-Euclidean lattices which should be the `right' objects associated to the loop group.} We refer to the main body of this paper (see in particular \S \ref{subIndHermitian}-\ref{subProHermitian}) for precise definitions and just offer a few comments on these notions here. The notion of ind-Euclidean lattices, or roughly, an increasing union of finite-rank lattices, is well-behaved with respect to taking limits of the corresponding finite-rank theta invariants. On the other hand, the theta-invariants for pro-Euclidean lattices require more care to define as a number of pathological phenomena arise with the~`naive' definition as a limit of the finite-rank quotients of the pro-Euclidean lattice. However, if the pro-Euclidean lattice satisfies a condition known as \emph{theta-finiteness} (see \S \ref{subThetaFiniteness}) one can in fact show that the naive definition of a theta invariant is in fact well-behaved and essentially the only one possible.

Finally, let us turn to the actual results of this paper. We focus on the parabolic arithmetic quotients \be{} \begin{array}{lllllll} \ovX^{\tau}_{\bb} & = & \wh{K} \setminus \gg_{\R}^{\tau} / \wh{\Gamma} \cap \bb & \; \text{ and } \; & \ovX^{\tau}_{\bb^-} & = & \wh{K} \setminus \gg_{\R}^{\tau} / \wh{\Gamma} \cap \bb^- \end{array} \ee where $\bb$ is the Kac--Moody analogue of a Borel subgroup (and $\bb^-$ is its opposite). When we have $\left \vert \tau \right \vert<1$, using Garland's reduction theory, one can show that the former space is the product of a finite-dimensional piece and a compact pro-unipotent piece. On the other hand, the latter space is really quite infinite-dimensional and our main theta-finiteness result only applies to a certain `polynomial' portion of this space.  Our main result is roughly the following (see Proposition \ref{prop:PsiKGamma} and Theorem \ref{thm:main-fin} for the precise version).

\begin{nthm}  To every $x \in \ovX^{\tau}_{\bb^-}$ we may attach a pro-Hermitian bundle on $\Spec(\zee)$. If we have $0 < \tau < 1$ and $x$ has a polynomial representative, then the corresponding pro-Hermitian bundle is also theta-finite. \end{nthm} 

In fact, to obtain a pro-Hermitian bundle one also needs to pick a representation of $\bb^-$ together with integral and Hermitian structures on this representation. In the main body of our paper, we work with irreducible highest-weight representations for which Garland (\cite{gar:loopAlgebras}) has constructed a natural integral form as well as a positive-definite Hermitian inner product. One can also consider the easier case of the adjoint representation using the same techniques as in this paper, the simplifications stemming from the fact that the dimensions of the finite rank quotients in the projective system corresponding to the adjoint representation grow linearly, and also that the lattice is quite explicitly given. In contrast, in the case of highest weight representations, the corresponding dimensions grow like the partition function and the lattices involved are quite subtle, containing contributions from the imaginary roots of the Kac--Moody root system and that are essentially described in terms of the Frenkel--Kac vertex operators or, equivalently, the `homogeneous' symmetric functions (see \S \ref{subsubIntegralFormZBasis}). 

To prove theta-finiteness, the polynomiality condition we impose requires $x$ to be chosen in $G(k[t, t^{-1}])$ rather than in $G(k((t)))$. Equivalently we require that it has an Iwasawa factorization with respect to $\bb^-.$
Our proofs rely on a few different ingredients: first, we need precise knowledge of how the finite-rank quotients in our projective system grow, which follows from certain estimates derived from the Weyl--Kac character formula and the theta-like behaviour of characters of highest weight modules. Next, we need to understand how the shortest length vectors in each of these pieces grows. To do this we rely on ideas essentially going back to H. Garland that appeared in his study of the convergence of (positive) Eisenstein series on loop groups (see \cite{gar:duke}). In particular, we need a linear-quadratic relation between the growth of the `central' and of the `classical' directions and a method of estimating norms of the action of unipotent elements in the highest-weight representation. Finally, we also use in a crucial way an estimate from the finite-dimensional world, namely a bound by Groenewegen (see \cite[Thm. 4.4]{gro:bound} or \cite[Prop. 2.6.2]{bost:book}) on the theta invariant of a Euclidean lattice in terms of its rank and of the shortest length of a non-zero vector. In a sense, our arguments circumvent much of the fine structure of the actual lattice.

A variant of our construction which will be described more fully in a future work is as follows. Attached to each element $x \in \ovX^{\tau}_{\bb}$ and some representation of $\bb$ we show how to construct a natural \emph{ind}-Hermitian bundle. In the case when the chosen representation is the adjoint representation of the Lie algebra of $\bb$ and assuming we have $0 < \tau <1$, one can then show the corresponding ind-Hermitian bundle has finite $h^0_{\vartheta}$. Similarly, if we choose $\tau > 1$, one may argue that the dual of the corresponding system, which is now a \emph{pro}-Hermitian bundle, is theta-finite. We do not at the moment know whether similar results also hold for a highest-weight representation. The issue is that when taking the dual lattice to Garland's integral form, we need better bounds than we currently know on the shortest vectors. Our difficulties stem exactly from trying to confront the intricate combinatorics of the imaginary root contributions head on.  

\subsection{Organization of the paper} In Section \ref{SecHermitianLattices}, we discuss the main features we need from Bost's theory of theta invariants for infinite-rank Euclidean lattices. While the main results of this paper concern pro-Hermitian bundles, we also included in \S \ref{subIndHermitian} a brief discussion of the simpler notion of ind-Hermitian bundles\footnote{Although it is not used in this paper, the duality between ind and pro-Hermitian bundles, discussed in \S \ref{SecHermitianLattices}, provides a construction, distinct from the one considered in \S \ref{SecProHermitian}, which we believe to be of interest.}. After reviewing these notions about theta-invariants, in \S \ref{SecLoopStuff} we give a fairly detailed description of the loop algebras and loop groups which we will be concerned with in this paper. In our discussion, we follow the work of Garland which emphasizes integral structures. Since this paper deals with aspects of Arakelov geometry and Kac--Moody theory, we give a presentation which tries to accommodate readers from both fields. For example, as our constructions of a pro-Hermitian bundle uses an integral structure on highest weight representations discovered by Garland in the late 70s, in \S  \ref{subsubIntegralFormZBasis} we give a bit more detail than we absolutely need about it. We thought the connection with symmetric functions sketched there might provide a concrete, if complicated (at least for us), way to approach the subject. Finally our main results are presented in \S \ref{SecProHermitian}.

\subsection{Acknowledgements} It is an honour to be able to dedicate this paper to the M. V. Subbarao on the occasion of the centennial of his birth. Not only was he instrumental in planting the seeds of the number theory community at the University of Alberta of which both authors are now a part, but the generosity of his family and their gift to the University of Alberta is the primary reason the authors were able to carry out their work together. Although it may seem far removed from his oeuvre, the partition function, one of Subbarao's favourite objects, plays a crucial role in this work. 

M.D. would like to thank Jean-Benoît Bost for fruitful conversations around this topic. M.P. was also partly supported by NSERC Discovery Grant RGPIN-2019-06112,  and he would also like to thank Howard Garland and Stephen D. Miller for a number of discussions around the topic of this paper and especially for their collaboration on Eisenstein series which implicitly informs this work. Finally we thank the referee for his or her careful reading of this paper.

\subsection{Notations} Let us make precise some of the notations used throughout this paper:

    \begin{itemize}
        \item $\mathbb{N}:= \{ 0, 1, 2, \ldots \} $ denotes the set of \textit{natural integers} $n \geqslant 0$, and $\mathbb{Z}$ the set of all integers;
        
        \item for any integers $a,b \in \mathbb{Z}$ with $a < b$, we denote by $\llbracket a, b \rrbracket$ the set of integers $k$ satisfying $a \leqslant k \leqslant b$;
        
        \item we denote by $A^{\times}$ the group of units of a ring $A$, but if $k$ is a field, we allow ourselves to write either $k^{\ast}$ or $k^{\times}$ for the group of non-zero elements;
        
        \item we denote by $\mathbb{R}_+^{\ast}:= \{ x \in \mathbb{R} \mid x > 0 \}$ the group of positive real numbers;
        
        \item the notation $\left< \cdot \right>$ means ``the subgroup generated by'' the elements in the brackets;
        
        \item the notation $\left< \cdot, \cdot \right>$ with two parameters will be used for bilinear forms and/or duality pairings.
    \end{itemize}

\section{Euclidean lattices in finite and infinite rank}
\label{SecHermitianLattices}

One of the central notions in Arakelov geometry is that of (finite-rank) Hermitian vector bundles over arithmetic curves, in particular over $\Spec \, \mathcal{O}_K$, where $K$ is a number field, and their connection to Euclidean lattices. In this section, we will review some aspects of this theory when  $K = \mathbb{Q}$, as this is the only case we consider in this paper. A treatment for general number fields will be considered in a future work. After presenting the finite-dimensional picture, we explain, following \cite{bost:book}, parts of Bost's extension of this theory to the infinite-dimensional setting. We refer to \emph{op. cit.} for more details.

    \subsection{\texorpdfstring{Finite-rank Hermitian vector bundles over $\Spec \, \mathbb{Z}$}{Finite-rank Hermitian vector bundles over Spec Z}}
    
        We begin this section by recalling some definitions and properties of Hermitian vector bundles over $\Spec \, \mathbb{Z}$.
        
        \spoint A \textit{Hermitian vector bundle} $\overline{E}$ over $\Spec \, \mathbb{Z}$ is the data 
        \be{}\label{dataHermitianVectBun} \begin{array}{ccc}
            \overline{E} & = & \left( E, \; \left \Vert \cdot \right \Vert \right)
        \end{array} \ee
        of a finitely generated, projective $\mathbb{Z}$-module $E$, with a Hermitian norm invariant by complex conjugation on the complex vector space $E_{\mathbb{C}} = E \otimes_{\mathbb{Z}} \mathbb{C}$. This definition coincides with that of a locally free sheaf over~$\Spec \, \mathbb{Z}$ with a Hermitian metric on the induced complex vector space, by taking the global sections. A survey of these objects is done in \cite[Ch. 1]{bost:book}, and in \cite{chen-moriwaki:book:arakelov}.
        
        \addtocounter{theorem}{1}
        
        \spoint Equivalently, the only infinite place of $\mathbb{Q}$ being real, a Hermitian vector bundle over $\Spec \, \mathbb{Z}$ can be seen as a \textit{Euclidean lattice}, meaning a free $\mathbb{Z}$-module of finite rank, together with a Euclidean norm on the real vector space $E_{\mathbb{R}} = E \otimes_{\mathbb{Z}} \mathbb{R}$.
        
        \addtocounter{theorem}{1}
        
        \spoint The \textit{rank} of a Hermitian vector bundle $\overline{E}$ is defined as the rank of $E$ as a $\mathbb{Z}$-module. Since $E$ is finitely generated, this is the same as the complex dimension of $E_{\mathbb{C}}$ or the real dimension of $E_{\mathbb{R}}$. We further say that $\overline{E}$ is a \textit{Hermitian line bundle} over $\Spec \, \mathbb{Z}$ if its rank equals $1$.
        
        \addtocounter{theorem}{1}
        
        \tpoint{Morphisms}\label{morphismsFiniteRankHermVB} Consider two Hermitian vector bundles $\overline{E}$ and $\overline{F}$ over $\Spec \, \mathbb{Z}$. A \textit{morphism} between~$\overline{E}$ and $\overline{F}$ is the datum of a $\mathbb{Z}$-linear map $\psi : E \longrightarrow F$. Note that the compatibility with the metrics comes from the fact that any $\mathbb{C}$-linear map in finite dimension is continuous. We denote the set of such morphisms by $\Hom_{\mathbb{Z}} \left( E, \, F \right)$. For any real number $\lambda > 0$, we further set
        \be{}\label{morphismBoundedNorm}
            \begin{array}{ccc}
                \Hom_{\mathbb{Z}}^{\leqslant \lambda} \left( \overline{E}, \, \overline{F} \right) & = & \left \{ \psi \in \Hom_{\mathbb{Z}} \left( E, \, F \right), \; \left \Vert \psi_{\mathbb{C}} \left( x \right) \right \Vert \; \leqslant \; \lambda \left \Vert x \right \Vert \text{ for any } x \in E_{\mathbb{C}} \right \}
            \end{array},
        \ee where $\psi_{\mathbb{C}} : E_{\mathbb{C}} \longrightarrow F_{\mathbb{C}}$ is the linear map induced by $f$ between the respective complexifications of $E$ and $F$. Such a morphism is said to be an \textit{isometry} if it satisfies $\left \Vert \psi_{\mathbb{C}} \left( x \right) \right \Vert = \left \Vert x \right \Vert$ for any $x \in E_{\mathbb{C}}$, in which case it is automatically bijective, as we are working with finite-rank vector bundles. We denote by $\vectZ$ the category of Hermitian vector bundles over $\Spec \, \zee$.
        
        \addtocounter{theorem}{1}

        \tpoint{Injective admissible morphisms} A morphism $\psi \in \Hom_{\mathbb{Z}}^{\leqslant 1} \left( \overline{E}, \, \overline{F} \right)$ is said to be \textit{injective admissible} if the underlying map $\psi : E \longrightarrow F$ is injective, with a torsion-free cokernel, and if it induces an isometry onto its image.
        
        \addtocounter{theorem}{1}

        \tpoint{Surjective admissible morphisms} A morphism $\psi \in \Hom_{\mathbb{Z}}^{\leqslant 1} \left( \overline{E}, \, \overline{F} \right)$ is said to be \textit{surjective admissible} it is surjective, and induces an isometry
        \be{}\label{isomSurjAdmMorph}
            \begin{array}{ccccc}
                \widetilde{\psi} & : & E/\ker \psi & \longrightarrow & F
            \end{array},
        \ee where the quotient space $E_{\mathbb{C}}/\left( \ker \psi \right)_{\mathbb{C}}$ is endowed with the quotient metric, defined by
        \be{}\label{defQuotientMetric}
            \begin{array}{ccc}
                \left \Vert x + \left( \ker \psi \right)_{\mathbb{C}} \right \Vert & = & \inf \left \{ \left \Vert x + y \right \Vert, \; y \in \left( \ker \psi \right)_{\mathbb{C}} \right \}
            \end{array}.
        \ee For more information on quotient metrics, the reader is referred to \cite[Sec. 1.1.3]{chen-moriwaki:book:arakelov}. Note that if $\psi: \ov{E} \longrightarrow \ov{F}$ is injective and admissible, then the transpose map $\psi^{\vee}: \overline{F}^{\vee} \longrightarrow \overline{E}^{\vee}$ is surjective and admissible.
        
        \addtocounter{theorem}{1}

        \tpoint{Covolume} Let $\overline{E}$ be a Hermitian vector bundle over $\Spec \, \mathbb{Z}$ of rank $r$, and $e_1$, \ldots, $e_r$ be a $\mathbb{Z}$-basis of $E$. Denote by $\left< \cdot, \cdot \right>$ the Euclidean inner product on $E_{\mathbb{R}}$. We define the \textit{covolume} of $\overline{E}$ by
        \be{}\label{defCovolume}
            \begin{array}{ccc}
                \covol \overline{E} & = & \sqrt{\det \left( \left( \left< e_i, \, e_j \right> \right)_{i,j=1,\ldots,r} \right)}
            \end{array}.
        \ee In other words, it is the square root of the determinant of the Gram matrix associated to the chosen basis. Note that this definition does not depend on the choice of $\mathbb{Z}$-basis of $E$, as two such bases are related by a integral matrix of determinant $1$.
        
        \addtocounter{theorem}{1}

        \tpoint{Arithmetic degree} There also exists a notion of degree of a line bundle, which takes into account the metric datum in a Hermitian vector bundle $\overline{E}$ over $\Spec \, \mathbb{Z}$. Let us first deal with the case of a Hermitian line bundle $\overline{L}$. Consider a non-zero element $s \in L$. Every element of $L$ can be uniquely written as $x s$, where $x$ is a rational, thereby giving an embedding $L \hookrightarrow \mathbb{Q}$, whose image is fractional ideal of $\mathbb{Q}$, and thus generated as a $\mathbb{Z}$-module by a rational number $y$. Denote by $n_p$ the $p$-adic valuation of $y$, which coincides with the opposite of the order of vanishing of $s$ at $p$. We then define the \textit{arithmetic degree} of $\overline{L}$ by
        \be{}\label{defArithmeticDegree}
            \begin{array}{ccc}
                \deghat \overline{L} & = & \sideset{}{}\sum\limits_{p \text{ prime}} \; n_p \; \log p \; - \; \log \left \Vert s \right \Vert
            \end{array}
        \ee where $\left \Vert \cdot \right \Vert$ is the Euclidean norm on $L_{\mathbb{R}}$. For a more general Hermitian vector bundle $\overline{E}$ of rank $r$, we set
        \be{}\label{defArithmeticDegreeHigherRank}
            \begin{array}{ccc}
                \deghat \overline{E} & = & \deghat \Lambda^{r} \; \overline{E}
            \end{array},
        \ee where $\Lambda^r \; \overline{E}$ denotes the top exterior power of $\overline{E}$, and is endowed with the determinant metric. This invariant is closely related to the covolume, as we have
        \be{}\label{relationCovolumeArithmeticDegree}
            \begin{array}{ccc}
                \deghat \overline{E} & = & - \log \covol \overline{E}
            \end{array}.
        \ee More on this and related notions can be found in \cite[Sec. 1.3]{bost:book}, and in \cite[Sec. 1]{freixas:article:arakelovNotes}. As a first example, one can take \cite[Sec. 1.1.3]{bost:book}, where, for any real number $\delta > 0$, the Hermitian line bundle $\overline{\mathcal{O}} \left( \delta \right)$ is defined by
        \be{}\label{defOdelta}
            \begin{array}{ccc}
                \overline{\mathcal{O}} \left( \delta \right) & = & \left( \mathbb{Z}, \; \left \Vert \cdot \right \Vert_{\delta} \right)
            \end{array},
        \ee the norm being characterized by $\left \Vert 1 \right \Vert_{\delta} = e^{- \delta}$. From \eqref{defArithmeticDegree}, we have $\deghat \overline{\mathcal{O}} \left( \delta \right) = \delta$.
        
        \addtocounter{theorem}{1}
        
        \tpoint{Dual lattice} For any Hermitian vector bundle $\overline{E}$, one can define the \textit{dual} bundle $\overline{E}^{\vee}$ by
        \be{}\label{defDualHermitianVB}
            \begin{array}{ccc}
                E^{\vee} & = & \Hom \left( E, \, \mathbb{Z} \right)
            \end{array},
        \ee whose complexification is endowed with the operator norm. It is a Hermitian vector bundle of the same rank as $\overline{E}$, and we can relate its arithmetic degree to that of $\overline{E}$, by
        \be{}\label{arithmeticDegreeDual}
            \begin{array}{ccc}
                \deghat \overline{E}^{\vee} & = & - \; \deghat \overline{E}
            \end{array}.
        \ee

    \subsection{Theta invariants in finite rank}
    
        The invariants we will be concerned with in this paper are called \textit{theta invariants}. Unlike the covolume and the arithmetic degree, we will see later that those can be generalized in an infinite-dimensional setting. Let us review them for finite-dimensional Hermitian vector bundles.
    
        \tpoint{\texorpdfstring{The invariants $h_{\vartheta}^0$ and $h_{\vartheta}^1$}{The invariants h0 and h1}}\label{thetaInvariantsFiniteRank} Let $\overline{E}$ be a Hermitian vector bundle over $\Spec \, \mathbb{Z}$. We set
        \be{}\label{defInvarianth0FiniteRank}
            \begin{array}{ccc}
                h_{\vartheta}^0 \left( \overline{E} \right) & = & \log \sideset{}{}\sum\limits_{v \in E} \; e^{- \pi \left \Vert v \right \Vert^2}
            \end{array},
        \ee the convergence of the sum above being guaranteed by the Poisson summation formula. Furthermore, note that this invariant also has an analytic interpretation, as the logarithm of the heat trace of an explicit multiple of the Laplacian on the quotient space $E_{\mathbb{R}}/E$. Using duality, we can now define the invariant $h_{\vartheta}^1$ as
        \be{}\label{defInvarianth1FiniteRank}
            \begin{array}{ccc}
                h_{\vartheta}^1 \left( \overline{E} \right) & = & h_{\vartheta}^0 ( \overline{E}^{\vee} )
            \end{array}.
        \ee
        
        \addtocounter{theorem}{1}
        
        \tpoint{The Riemann--Roch formula} As a further consequence of the Poisson summation formula, we can compare the invariants $h_{\vartheta}^0$ and $h_{\vartheta}^1$ to the covolume, as
        \be{}\label{RRFormulaFiniteRank}
            \begin{array}{ccc}
                h_{\vartheta}^0 \left( \overline{E} \right) \; - \; h_{\vartheta}^1 \left( \overline{E} \right) & = & \deghat \overline{E}
            \end{array}.
        \ee
        
        \addtocounter{theorem}{1}
        
        \tpoint{\texorpdfstring{Upper-bound on $h_{\vartheta}^0$}{Upper-bound on h0}} Let $\overline{E}= (E, || \cdot ||)$ be a Hermitian vector bundle over $\Spec \, \mathbb{Z}$. We denote by \be{ }\label{defn:shortestvector} \lambda_1 \text { or more precisely } \lambda_1(\overline{E}) \; = \; \lambda_1(E, || \cdot ||)  \ee  the first minimum of $\overline{E}$, meaning the shortest norm of a non-zero vector. To emphasize the norm chosen, we abuse notation slightly and write $\lambda_1(\overline{E}, || \cdot ||)$ for this quantity as well. 
        
        \begin{theorem} \label{thmGroenewegen} \cite[Prop. 2.6.2]{bost:book}, \cite[Prop. 4.4]{gro:bound}
            Denote by $r$ the rank of $E$. We have \be{}\label{upperBoundh0} \begin{array}{lllll} h_{\vartheta}^0 \left( \overline{E} \right) & \leqslant & C \left( r, \, \lambda_1 \right) & = & \displaystyle 3^r \left( \pi \lambda_1^2 \right)^{-r/2} \int_{\pi \lambda_1^2}^{+ \infty} \; u^{r/2} \, e^{-u} \; \mathrm{d}u \end{array}. \ee Assuming we have $\lambda_1 > \left( r / 2\pi \right)^{1/2}$, we further have \be{}\label{inequalityConstantCrlambda} \begin{array}{lll} C \left( r, \, \lambda_1 \right) & \leqslant & \displaystyle 3^r \left( 1 - \frac{r}{2 \pi \lambda_1^2} \right)^{-1} e^{- \pi \lambda_1^2} \end{array}. \ee
        \end{theorem}

    \subsection{\texorpdfstring{Ind-Hermitian vector bundles over $\Spec \, \mathbb{Z}$}{Ind-Hermitian vector bundles over Spec Z}}
    \label{subIndHermitian}
        The first generalization of the notion of Hermitian vector bundles to an infinite-dimensional setting is that of an ind-Hermitian vector bundle. The reader will find a deeper study of these in \cite[Sec. 5]{bost:book}.
        
        \spoint An \textit{Ind-Hermitian} vector bundle $\overline{E}$ over $\Spec \, \mathbb{Z}$ is the data
        \be{}\label{dataHermitianVectBun} 
            \begin{array}{ccc}
                \overline{E} & = & \left( E, \; \left \Vert \cdot \right \Vert \right)
            \end{array},
        \ee of a countably generated, projective $\mathbb{Z}$-module $E$, and of a pre-Hilbertian norm $\left \Vert \cdot \right \Vert$ invariant by complex conjugation on $E_{\mathbb{C}} = E \otimes_{\mathbb{Z}} \mathbb{C}$.

        \tpoint{Morphisms of ind-Hermitian vector bundles} A \textit{morphism} of ind-Hermitian vector bundles is the datum of a $\mathbb{Z}$-linear map $\psi : E \longrightarrow F$ which induces a continuous complex linear map $\psi_{\mathbb{C}} : E_{\mathbb{C}} \longrightarrow F_{\mathbb{C}}$ on the complexifications. Note that, unlike the finite rank case, continuity is no longer automatic here. We denote the set of such morphisms by $\Hom \left( \overline{E}, \overline{F} \right)$ and we further set, for any $\lambda >0$,
        \be{}\label{morphIndHermBoundedNorm}
            \begin{array}{ccc}
                \Hom_{\mathbb{Z}}^{\leqslant \lambda} \left( \overline{E}, \, \overline{F} \right) & = & \left \{ \psi \in \Hom_{\mathbb{Z}} \left( E, \, F \right), \; \left \Vert \psi_{\mathbb{C}} \left( x \right) \right \Vert \; \leqslant \; \lambda \left \Vert x \right \Vert \text{ for any } x \in E_{\mathbb{C}} \right \}
            \end{array}.
        \ee An element $\psi \in \Hom \left( \overline{E}, \overline{F} \right)$ is called an \textit{isometry of ind-Hermitian vector bundles} if it is bijective and if we have $\left \Vert \psi_{\mathbb{C}} \left( x \right) \right \Vert = \left \Vert x \right \Vert$ for any $x \in E_{\mathbb{C}}$. We denote by $\indvectZ$ the category of ind-Hermitian vector bundles over $\Spec \, \mathbb{Z}$.

        \tpoint{Admissible inductive systems} An \textit{admissible inductive system} of Hermitian vector bundles is a complex of (finite rank) Hermitian vector bundles
        \be{}\label{inductiveSystem} 
            \begin{array}{ccccccccccccc}
                \overline{E}_{\bullet} & : & \overline{E}_0 & \overset{j_0}{\longhookrightarrow} & \overline{E}_1 & \overset{j_1}{\longhookrightarrow} & \ldots &  \overset{j_{n-1}}{\longhookrightarrow} & \overline{E}_n & \overset{j_n}{\longhookrightarrow} & \overline{E}_{n+1} & \ldots
            \end{array}
        \ee over $\Spec \, \mathbb{Z}$, where each morphism $j_n$ is injective admissible.

        \spoint To an inductive system as \eqref{inductiveSystem}, we can attach an ind-Hermitian vector bundle, by setting
        \be{}\label{indHermitianFromInductiveSystem}
            \begin{array}{ccccc}
                E & = & \lim\limits_{\underset{n}{\longrightarrow}} \; E_n & = & \sideset{}{}\bigsqcup\limits_n \; E_n \; \; / \; \;  \sim
            \end{array},
        \ee where the equivalence relation $\sim$ on the disjoint union is defined by $x_k \sim x_n$ for any $x_k \in E_k$ and $x_n \in E_n$ if both elements are equal in some larger $E_m$ after composition by the maps $j$. The complexification $E_{\mathbb{C}}$ is endowed with the unique norm $\left \Vert \cdot \right \Vert$ such that the maps $E_{n, \mathbb{C}} \hookrightarrow E_{\mathbb{C}}$ are injective admissible. Conversely, any ind-Hermitian vector bundle comes from an (and possibly many different) inductive systems.

        \spoint A morphism $\varphi : \overline{E}_{\bullet} \longrightarrow \overline{F}_{\bullet}$ of inductive systems is the datum of a collection $\left( \varphi_n \right)_n$ of morphisms making the following diagram commutative
        \be{}\label{morphismInductiveSystem} 
            \begin{array}{lllllllllllll}
                \overline{E}_{\bullet} & : & \overline{E}_0 & \overset{j_0}{\longhookrightarrow} & \overline{E}_1 & \overset{j_1}{\longhookrightarrow} & \ldots &  \overset{j_{n-1}}{\longhookrightarrow} & \overline{E}_n & \overset{j_n}{\longhookrightarrow} & \overline{E}_{n+1} & \longhookrightarrow & \ldots \\[0.5em]
                \downarrow \varphi & & \downarrow \varphi_0 & & \downarrow \varphi_1 & & & & \downarrow \varphi_n & & \downarrow \varphi_{n+1} \\
                \overline{F}_{\bullet} & : & \overline{F}_0 & \overset{j_0'}{\longhookrightarrow} & \overline{F}_1 & \overset{j_1'}{\longhookrightarrow} & \ldots &  \overset{j_{n-1}'}{\longhookrightarrow} & \overline{F}_n & \overset{j_n'}{\longhookrightarrow} & \overline{F}_{n+1} & \longhookrightarrow & \ldots
            \end{array}
        \ee Such a morphism induces a morphism of the associated ind-Hermitian vector bundles as in \eqref{indHermitianFromInductiveSystem}.

    \subsection{\texorpdfstring{Pro-Hermitian vector bundles over $\Spec \, \mathbb{Z}$}{Pro-Hermitian vector bundles over Spec Z}}
     \label{subProHermitian}    
        The more complicated notion of Hermitian vector bundles is that of \textit{pro-Hermitian vector bundle}. As the name suggests, there is a close link to projective systems. However, Bost offers in \cite[Sec. 5.1.2, 5.1.4]{bost:book} definitions which do not mention projective systems explicitly. In \cite[Sec. 5.4.2]{bost:book}, morphisms between pro-Hermitian vector bundles are defined, using the aforementioned definitions. The category of pro-Hermitian vector bundles over $\Spec \, \mathbb{Z}$ is denoted by $\provectZ$. To lighten the presentation, we will only consider the incarnation of pro-Hermitian vector bundles from \textit{admissible projective systems}.
        
        \tpoint{Admissible projective systems}\label{admissibleProjSyst}  An \textit{admissible projective system} of Hermitian vector bundles is a complex of (finite-rank) Hermitian vector bundles
        \be{}\label{projectiveSystem} 
            \begin{array}{ccccccccccccc}
                \overline{E}_{\bullet} & : & \overline{E}_0 & \overset{q_0}{\longtwoheadleftarrow} & \overline{E}_1 & \overset{q_1}{\longtwoheadleftarrow} & \ldots &  \overset{q_{n-1}}{\longtwoheadleftarrow} & \overline{E}_n & \overset{q_n}{\longtwoheadleftarrow} & \overline{E}_{n+1} & \ldots
            \end{array}
        \ee over $\Spec \, \mathbb{Z}$, where each morphism $q_n$ is surjective admissible.

        \spoint\label{proHermitianFromProjSyst} To a projective system as \eqref{projectiveSystem}, we can attach a \textit{pro-Hermitian vector bundle}, by setting
        \be{}\label{proHermitianFromProjectiveSystem}
            \begin{array}{ccccc}
                E & = & \lim\limits_{\underset{n}{\longleftarrow}} \; E_n & = & \left \{ \left( x_n \right)_n \in \sideset{}{}\prod\limits_{n \geqslant 0} \; E_n, \; q_k \left( x_{k+1} \right) = x_k \text{ for every } k \in \mathbb{N} \right \}
            \end{array},
        \ee whose complexification is endowed with a norm described in \cite[Sec. 5.3.2]{bost:book}, which we omit here to shorten this survey. A difference between ind- and pro-Hermitian vector bundles is that \eqref{proHermitianFromProjectiveSystem} does not in general lead to a countably generated $\mathbb{Z}$-module.

        \tpoint{Morphisms of projective systems} A collection $\left( \psi_n \right)_n$ of morphisms such that the diagram
        \be{}\label{morphismProjectiveSystem} 
            \begin{array}{lllllllllllll}
                \overline{E}_{\bullet} & : & \overline{E}_0 & \overset{q_0}{\longtwoheadleftarrow} & \overline{E}_1 & \overset{q_1}{\longtwoheadleftarrow} & \ldots &  \overset{q_{n-1}}{\longtwoheadleftarrow} & \overline{E}_n & \overset{q_n}{\longtwoheadleftarrow} & \overline{E}_{n+1} & \longtwoheadleftarrow & \ldots \\[0.5em]
                & & \downarrow \psi_0 & & \downarrow \psi_1 & & & & \downarrow \psi_n & & \downarrow \psi_{n+1} \\
                \overline{F}_{\bullet} & : & \overline{F}_0 & \overset{q_0'}{\longtwoheadleftarrow} & \overline{F}_1 & \overset{q_1'}{\longtwoheadleftarrow} & \ldots &  \overset{q_{n-1}'}{\longtwoheadleftarrow} & \overline{F}_n & \overset{q_n'}{\longtwoheadleftarrow} & \overline{F}_{n+1} & \longtwoheadleftarrow & \ldots
            \end{array}
        \ee is commutative induces a morphism $\psi : \overline{E} \longrightarrow \overline{F}$ of the associated pro-Hermitian vector bundles \eqref{proHermitianFromProjectiveSystem}, in the sense of \cite[Sec. 5.4.2]{bost:book}. Such a morphism is an \textit{isometry} in the sense of \cite[Sec. 5.1.2]{bost:book} if, for $n$ large enough, every $\psi_n$ is an isometry, as in \ref{morphismsFiniteRankHermVB}. Note that we include bijectivity in the term ``isometry''.

    \subsection{Duality} As explained in \cite[Sec. 5.5]{bost:book}, there is a correspondence between ind- and pro-Hermitian vector bundles over $\Spec \, \mathbb{Z}$, using duality.
    
        \spoint Consider an admissible inductive system $\overline{E}_{\bullet}$ as in \eqref{inductiveSystem}. One can then consider the system
        \be{}\label{projectiveDualSystem} 
            \begin{array}{ccccccccccccc}
                \overline{E}_{\bullet}^{\vee} & : & \overline{E}_0^{\vee} & \overset{q_0}{\longtwoheadleftarrow} & \overline{E}_1^{\vee} & \overset{q_1}{\longtwoheadleftarrow} & \ldots &  \overset{q_{n-1}}{\longtwoheadleftarrow} & \overline{E}_n^{\vee} & \overset{q_n}{\longtwoheadleftarrow} & \overline{E}_{n+1}^{\vee} & \ldots
            \end{array},
        \ee of finite-rank Hermitian vector bundles given by the dual of each $\overline{E}_n$, where the map $q_n : E_{n+1}^{\vee} \longrightarrow E_n^{\vee}$ is the transpose of $j_n : E_n \longrightarrow E_{n+1}$. It is then explained in \cite[Sec. 5.5.1]{bost:book} that the system \eqref{projectiveDualSystem} is projective admissible. Furthermore, taking the transpose of a morphism (resp. an isometry) of inductive systems yields a morphism (resp. an isometry) of projective systems.
        
        \spoint Consider an admissible projective system $\overline{F}_{\bullet}$ as in \eqref{projectiveSystem}. One can then consider the system
        \be{}\label{inductiveDualSystem}
            \begin{array}{ccccccccccccc}
                \overline{F}_{\bullet}^{\vee} & : & \overline{F}_0^{\vee} & \overset{j_0}{\longhookrightarrow} & \overline{F}_1^{\vee} & \overset{j_1}{\longhookrightarrow} & \ldots &  \overset{j_{n-1}}{\longhookrightarrow} & \overline{F}_n^{\vee} & \overset{j_n}{\longhookrightarrow} & \overline{F}_{n+1}^{\vee} & \ldots
            \end{array},
        \ee of finite rank Hermitian vector bundles given by the dual of each $\overline{F}_n$, where the map $j_n : F_n^{\vee} \longrightarrow F_{n+1}^{\vee}$ is the transpose of $q_n : F_{n+1} \longrightarrow F_n$. Using \cite[Sec. 5.5.1]{bost:book}, the system \eqref{inductiveDualSystem} is inductive admissible. Furthermore, taking the transpose of a morphism (resp. an isometry) of projective systems yields a morphism (resp. an isometry) of inductive systems.

    \subsection{Theta-finiteness and theta invariants in infinite rank}
    \label{subThetaFiniteness}
        We will now see how to extend the definition of the theta invariants $h_{\vartheta}^0$ and $h_{\vartheta}^1$ from \ref{thetaInvariantsFiniteRank} to ind- and pro-Hermitian vector bundles.
        
        \tpoint{\texorpdfstring{Invariant $h_{\vartheta}^0$ of an ind-Hermitian vector bundle}{Invariant h_theta0 of an ind-Hermitian vector bundle}} The first and most direct generalization of theta invariants to the infinite-dimensional setting has to do with ind-Hermitian vector bundles. Consider an ind-Hermitian vector bundle $\overline{E}$ obtained from an admissible inductive system 
        \be{}
            \begin{array}{ccccccccccccc}
                \overline{E}_{\bullet} & : & \overline{E}_0 & \overset{j_0}{\longhookrightarrow} & \overline{E}_1 & \overset{j_1}{\longhookrightarrow} & \ldots &  \overset{j_{n-1}}{\longhookrightarrow} & \overline{E}_n & \overset{j_n}{\longhookrightarrow} & \overline{E}_{n+1} & \ldots
            \end{array}
        \ee of finite-rank Hermitian vector bundles over $\Spec \, \mathbb{Z}$. We get an increasing sequence
        \be{}\label{increasingSequenceIndHermVB}
            \begin{array}{ccccccccccc}
                0 & \leqslant & h_{\vartheta}^0 \left( \overline{E}_0 \right) & \leqslant & h_{\vartheta}^0 \left( \overline{E}_1 \right) & \leqslant & \ldots & \leqslant & h_{\vartheta}^0 \left( \overline{E}_n \right) & \leqslant & \ldots
            \end{array},
        \ee which either converges, or diverges to $+ \infty$. We then set
        \be{}\label{h0thetaIndHermitian}
            \begin{array}{ccccc}
                h_{\vartheta}^0 \left( \overline{E} \right) & = & \lim\limits_{n \rightarrow + \infty} h_{\vartheta}^0 \left( \overline{E}_n \right) & \in & \left[ 0, \, + \infty \right]
            \end{array}.
        \ee This limit can be expressed as the following sum
        \be{}\label{h0thetaIndHermitianSeries}
            \begin{array}{ccccc}
                h_{\vartheta}^0 \left( \overline{E} \right) & = & \log \; \sideset{}{}\sum\limits_{v \in E} \; e^{- \pi \left \Vert v \right \Vert^2}
            \end{array}.
        \ee Using \cite[Prop. 6.1.1]{bost:book} and the reflexivity property \cite[Prop. 5.5.1]{bost:book}, one may define the $h_{\vartheta}^1$ of a pro-Hermitian vector bundle as the $h_{\vartheta}^0$ of its dual. Note that there would be a twist by the line bundle $\omega_{\mathcal{O}_K/\mathbb{Z}}$ for general number fields.
        
        %By duality, one may define the $h_{\vartheta}^1$ invariant of a Pro-Hermitian vector bundle as the $h_{\vartheta}^0$ of its dual. For more general number fields, a twist by a line bundle is needed here, this being explained in \cite[Prop. 6.1.1]{bost:book}.

        \tpoint{Theta invariants of a pro-Hermitian vector bundle} Defining the invariant $h_{\vartheta}^0$ of a pro-Hermitian vector bundle requires more care. Let us briefly see why. Considering such a bundle $\overline{F}$, coming from an admissible projective system
        \be{}\label{projectiveSystemTheta} 
            \begin{array}{ccccccccccccc}
                \overline{F}_{\bullet} & : & \overline{F}_0 & \overset{q_0}{\longtwoheadleftarrow} & \overline{F}_1 & \overset{q_1}{\longtwoheadleftarrow} & \ldots &  \overset{q_{n-1}}{\longtwoheadleftarrow} & \overline{F}_n & \overset{q_n}{\longtwoheadleftarrow} & \overline{F}_{n+1} & \ldots
            \end{array},
        \ee the sequence $\left( h_{\vartheta}^0 \left( F_n \right) \right)_n$ does not necessarily converge. Furthermore, the convergence and possible limit of this sequence could depend on the projective system, even if they yield isometric pro-Hermitian vector bundles, in the sense of \cite[Sec. 5.1.2]{bost:book}. Additionally, the tentative definition
        \be{}\label{tentativeDefh0thetaPro}
            \begin{array}{ccc}
                \text{`` } h_{\vartheta}^0 \left( \overline{F} \right) & = & \sideset{}{}\sum\limits_{v \in F} \; e^{- \pi \left \Vert v \right \Vert^2} \text{ ''}
            \end{array}
        \ee cannot actually be made, as $F$ is in general uncountable. In \cite[Sec. 6.2.2]{bost:book}, Bost defines instead two invariants: a lower theta invariant $\underline{h}_{\vartheta}^0 \left( \overline{F} \right)$, and an upper one $\overline{h}_{\vartheta}^0 \left( \overline{F} \right)$, which satisfy
        \be{}\label{comparisonLowerUpperThetaPro}
            \begin{array}{ccc}
                \underline{h}_{\vartheta}^0 \left( \overline{F} \right) & \leqslant & \overline{h}_{\vartheta}^0 \left( \overline{F} \right)
            \end{array}.
        \ee Note that these are defined using the characterization of pro-Hermitian vector bundles without explicit projective systems, as in \cite[Sec. 5.1.2]{bost:book}. A particular case of interest will be when the two upper and lower invariants coincide.

        \tpoint{Strong summability condition} Consider a projective system as in \eqref{projectiveSystemTheta}. For any integer $n \geqslant 0$, we consider the finite-rank Hermitian vector bundle
        \be{}\label{kernelProjectiveSystem}
            \begin{array}{ccc}
                \ker q_n & \subset & F_{n+1}
            \end{array}
        \ee whose complexification is endowed with the restriction of the metric of $\overline{F}_{n+1}$. The resulting Hermitian vector bundle is denoted by $\overline{\ker q_n}$. The projective system \eqref{projectiveSystemTheta} is \textit{summable} if we have
        \be{}\label{summabilityCondition}
            \begin{array}{ccccc}
                \text{\textbf{Sum}} \left( \overline{F}_{\bullet} \right) & : & \sideset{}{}\sum\limits_{n \geqslant 0} \; h_{\vartheta}^0 \left( \overline{\ker q_n} \right) & < & + \infty.
            \end{array}
        \ee As explained in \cite[Thm. 7.3.4]{bost:book}, this condition is enough to guarantee that we have
        \be{}\label{upperEqualsLimit}
            \begin{array}{ccc}
                \overline{h}_{\vartheta}^0 \left( \overline{F} \right) & = & \lim\limits_{n \rightarrow + \infty} \; h_{\vartheta}^0 \left( \overline{F}_n \right)
            \end{array}
        \ee but it does not imply equality with the lower theta invariant. To get that property, a more demanding condition is required. The projective system \eqref{projectiveSystemTheta} is said to be \textit{strongly summable} if we have
        \be{}\label{strongSummabilityCondition}
            \begin{array}{ccccc}
                \text{\textbf{Sum}} \left( \overline{F}_{\bullet} \otimes \overline{\mathcal{O}} \left( \varepsilon \right) \right) & : & \sideset{}{}\sum\limits_{n \geqslant 0} \; h_{\vartheta}^0 \left( \overline{\ker q_n} \otimes \overline{\mathcal{O}} \left( \varepsilon \right) \right) & < & + \infty.
            \end{array}
        \ee for some real number $\varepsilon > 0$. In that case, inequality \eqref{comparisonLowerUpperThetaPro} becomes an equality, and the common value of the lower and upper theta invariants is also given by \eqref{upperEqualsLimit}. We then set
        \be{}\label{h0thetaProHermVB}
            \begin{array}{ccccccc}
                h_{\vartheta}^0 \left( \overline{F} \right) & = & \underline{h}_{\vartheta}^0 \left( \overline{F} \right) & = & \overline{h}_{\vartheta}^0 \left( \overline{F} \right) & = & \lim\limits_{n \rightarrow + \infty} \; h_{\vartheta}^0 \left( \overline{F}_n \right).
            \end{array}
        \ee Note that the lower and upper theta invariants only depend on the isometry class of the pro-Hermitian vector bundle, not on the projective system yielding it. Therefore, strong summability related to any projective system allows the definition of the $h_{\vartheta}^0$ invariant by the common value in \eqref{comparisonLowerUpperThetaPro}. The fact that this invariant can be computed as a limit of finite-rank theta invariants, seen as the last part of \eqref{h0thetaProHermVB}, is \textit{a priori} only verified for strongly summable projective systems.

        \tpoint{Theta-finite Pro-Hermitian vector bundles} A pro-Hermitian vector bundle $\overline{F}$ over $\Spec \, \mathbb{Z}$ is said to be \textit{theta-finite} if the condition
        \be{}\label{thetaFiniteness}
            \begin{array}{ccccc}
                \text{\textbf{Sum}} \left( \overline{F}_{\bullet} \otimes \overline{\mathcal{O}} \left( - \frac{1}{2} \log t \right) \right) & : & \sideset{}{}\sum\limits_{n \geqslant 0} \; h_{\vartheta}^0 \left( \overline{\ker q_n} \otimes \overline{\mathcal{O}} \left( - \frac{1}{2} \log t \right) \right) & < & + \infty.
            \end{array}
        \ee is satisfied for every real number $t > 0$. In that case, we can define the \textit{theta function} of $\overline{F}$ by
        \be{}\label{defThetaFunctionProHermVB}
            \begin{array}{ccccc}
                \vartheta_{\overline{F}} \left( t \right) & = & h_{\vartheta}^0 \left( \overline{F} \otimes \overline{\mathcal{O}} \left( - \frac{1}{2} \log t \right) \right) & = & \lim\limits_{n \rightarrow + \infty} \; \sideset{}{}\sum\limits_{v \in F_n} \; e^{- \pi t \left \Vert v \right \Vert_{F_n}^2}
            \end{array}
        \ee for every $t>0$, where the projective system $\overline{F}_{\bullet}$ has $\overline{F}$ as projective limit.

\section{Loop algebras and loop groups (\emph{alias} affine Kac--Moody algebras and groups) } \label{SecLoopStuff}

In this section, after first fixing our notation for finite-dimensional Lie algebras in \S \ref{subFiniteDimensionalAlgebras}, we recall two constructions for the main infinite-dimensional Lie algebra of interest to us, the loop algebra, or more precisely a Kac--Moody algebra of untwisted affine type. After introducing the basic features of these algebras in \S \ref{subSecLoopAlgebraUniversalCovering}, we focus on their highest weight representations in \S \ref{subGarlandRep}, and describe their integral and Hermitian structures in \S \ref{subIntegralFormGarlandRep}. After these preliminaries, we introduce the loop groups relevant to our work, following Garland's construction (\textit{cf}. \cite{gar:loopGroups}). To orient the expert, let us note our groups are the `maximal' or `complete' loop groups, where the completion is in the direction of the positive Borel, \textit{i.e.} in $t$ and not in~$t^{-1}$, where $t$ is the loop parameter, and that our groups are further twisted by the loop rotation (whose action is by `contracting' high, positive powers of the loop).

\subsection{Finite-dimensional Lie algebras}  \label{subFiniteDimensionalAlgebras}
Let us begin with a quick review of the theory of finite-dimensional Lie algebras, as they constitute the first step towards defining loop algebras.

\tpoint{Root systems} \label{subsubFinRts} Let $\mf{g}$ be a complex, simple Lie algebra with a Cartan subalgebra $\mf{h}$ of rank $\ell$. Decomposing $\mf{g}$ with respect to the adjoint action of $\mf{h}$, we obtain $\mf{g} = \mf{h} \oplus \bigoplus_{\alpha} \mf{g}^{\alpha}$, with $\alpha \in \mf{h}^*$ and \be{wt-space} \begin{array}{lll} \mf{g}^{\alpha} & = & \{ X \in \mf{g} \mid [H, X] = \alpha(H) X \text{ for all } H \in \mf{h} \}. \end{array} \ee An element $\alpha$ as above is said to be a \textit{root} of $\mf{g}$ if we have $\mf{g}^{\alpha} \neq 0$. In this case, the space $\mf{g}^{\alpha}$ is one-dimensional. We denote by $\rts \subset \mathfrak{h}^*$ the set of roots, and the pair $\left( \mathfrak{h}^*, \rts \right)$ satisfy the axioms of an irreducible root system (\textit{cf}. \cite[Ch. VI, \S 1]{bourbaki:456}). Recall that such a root system has a basis $\Pi:= \{ \alpha_1, \ldots, \alpha_{\ell} \}$ of simple roots and that every $\alpha \in \rts$ is an integral linear combination $\alpha = \sum_{i =1}^{\ell} m_i \alpha_i$, where we either have $m_i \geqslant 0$ for every $i$, or $m_i \leqslant 0$ for every $i$.  In the former case, we say the root is positive and in the latter that it is negative. We also write $\rts_{\pm}$ for the set of positive/negative roots; the \emph{height} of any positive root $\alpha \in \rts_+$ is defined as the positive integer $\mathrm{ht}(\alpha) = \sum_{i=1}^{\ell} m_i$. There exists a unique root (\textit{cf}. \cite[Ch. VI, \S 8]{bourbaki:456}) of maximal height $\alpha_0 \in \rts_+$, called the \textit{highest root}, which we write as \be{highest-root} \begin{array}{lll} \alpha_0 & = & \sideset{}{}\sum\limits_{i=1}^{\ell} d_i \, \alpha_i. \end{array} \ee We have $d_i >0$ for every integer $i=1, \ldots, \ell$. Let us denote by $\mathscr{Q}$ the $\zee$-lattice spanned by $\Pi$ and by $\mathscr{Q}_+$ the set of $\zee_{\geq0}$-linear combinations of $\Pi$.

\tpoint{Normalized Killing form and coroots} \label{subsubKillingCoroots}Denote by $(\cdot, \cdot)_{\kil} : \mf{g} \times \mf{g} \rr \C$ the Killing form. Restricted to $\mf{h}$, it is a non-degenerate, symmetric, bilinear form. Hence, with respect to it, we may define the dual element $H'_{\alpha} \in \mf{h}$ of any root $\alpha \in \rts$. For any $\alpha, \beta \in \mf{h}^*$ we set $(\alpha, \beta)_{\kil}:= (H'_{\alpha}, H'_{\beta})_{\kil} $. When considering affine Lie algebras, we will also need the renormalization of the Killing form on $\mf{h}$ (and also $\mf{h}^*$) as \be{norm:kill} \begin{array}{lll} (\cdot , \cdot) & := & \displaystyle \frac{2}{(\alpha_0, \alpha_0)_{\kil}} \; (\cdot, \cdot)_{\kil}. \end{array} \ee For any $\alpha \in \rts$, we define the corresponding \emph{coroot} as  \be{coroot} \begin{array}{lllll} \alpha^{\vee} & := & \displaystyle \frac{2}{(\alpha, \alpha)_{\kil}} \; H'_{\alpha} & \in & \mf{h}, \end{array} \ee and let $\rts^{\vee}:= \{ \alpha^{\vee} \mid \alpha \in \rts \}$ be the set of coroots. We pick a basis $\Pi^{\vee}:= \{ \alphav_1, \ldots, \alphav_{\ell} \}$ of $\rts^{\vee}$, with respect to which we define the sets of positive and negative coroots $\crts_+$ and $\crts_-$, respectively. The pair $\left( \rts^{\vee}, \mf{h} \right)$ also forms a root system.

\tpoint{Cartan matrix} \label{subsubCartanMatrix} Let $\la \cdot, \cdot \ra$ denote the natural pairing $\mf{h} \times \mf{h}^* \rr \C.$ Define the Cartan matrix of $\mf{g}$ as  \be{cartan-matrix} \begin{array}{lllllllllll} A & := & (A_{ij}) & \text{ where } & A_{ij} & = & \la \alphav_i, \alpha_j \ra & = & 2 \; \frac{(\alpha_j, \alpha_i)_{\kil}}{(\alpha_i, \alpha_i)_{\kil}} & = &  2 \; \frac{(\alpha_j, \alpha_i)}{(\alpha_i, \alpha_i)} . \end{array} \ee

\tpoint{Weyl group} \label{subsubWeylGroup} We write $W$ for the Weyl group of the root system $\rts$ (or $\crts$). It is a finite group generated by the simple reflections $s_i$ which act on $\mf{h}^*$ as $s_i(\lambda) = \lambda - \la \lambda, \alphav_i \ra \alpha_i.$ In fact, $(W, S)$ is a Coxeter system in the sense of \cite[Ch. IV]{bourbaki:456}. Denote the length function on this Coxeter group as $\ell: W \rr \mathbb{N}\cup\{ 0 \}.$

\tpoint{Chevalley form} \label{subsubChevalleyBasisIntegralForm}
In what follows, we need to pick a specific basis of $\mf{g}$ with respect to which the structure coefficients are integral. As explained in \cite[Ch. 1, Thm. 1]{steinberg:book:lecturesChevalleyGroups}, for every $\alpha \in \rts$, we may pick an element $E_{\alpha} \in \mathfrak{g}^{\alpha}$ such that we have               
    \[ \begin{array}{rccccrcccl}
            1. & \left[ E_{\alpha}, E_{- \alpha} \right] & = & \alphav  \hspace{0.3em} & & 2. & \left[ E_{\alpha}, E_{\beta} \right] & = & \pm \left( r+1 \right) E_{\alpha + \beta} & \text{if } \alpha + \beta \in \Delta \text{ and $r$ is the largest } \\
            &&&&&&&&& \text{ integer $k$ with $a-k\beta \in \rts(A)$} \\[0.5em]
                        
            3. & \left[ \alpha^{\vee}, E_{\beta} \right] & = & 2 \frac{\left( \alpha, \beta \right)}{\left( \alpha, \alpha \right)} E_{\beta} & \hspace{0.3em} & 4. & \left[ E_{\alpha}, E_{\beta} \right] & = & 0 & \text{if $2$ is not satisfied}
        \end{array} \]

\noindent The set \be{ChevBasis} \mathcal{B}_{\mathrm{Ch}} ( \mathfrak{g} ) = \left \{ E_{\alpha}, \alpha \in \rts( A ) \right \} \cup \left \{ \alphav_1, \dots, \alphav_{\ell} \right \} \ee is called a \emph{Chevalley basis} of $\mathfrak{g}$ and enjoys a number of nice properties, including the integrality mentioned above. We refer to \cite[Thm. 1, p.7]{steinberg:book:lecturesChevalleyGroups} for more details. Using this basis $\mathcal{B}_{\mathrm{Ch}} ( \mathfrak{g} )$, one can define a $\zee$-form of the enveloping algebra of $\mathfrak{g}$ and consruct a $\zee$-basis of it explicitly (\cite[Ch. 2]{steinberg:book:lecturesChevalleyGroups}. As we need these constructions for the loop algebra, we review them in that context later (see \S \ref{subIntegralFormGarlandRep}).

\subsection{Loop algebras and (untwisted) affine Kac--Moody algebras }
\label{subSecLoopAlgebraUniversalCovering}
In this section, we give two constructions of the main infinite-dimensional Lie algebra under study here. The first is often called the loop presentation and the second the Kac--Moody presentation. Our approach follows \cite{gar:loopGroups, gar:loopAlgebras} and \cite{kac:book}.
          
\newcommand{\Lg}{L\mf{g}}
\newcommand{\gext}{\widehat{\mf{g}}}

\tpoint{Loop presentation} \label{subsubLoop} The \textit{loop algebra} $\widetilde{\mathfrak{g}}$ of $\mathfrak{g}$ is defined as a $k$-vector space by $\widetilde{\mathfrak{g}} = \mathfrak{g} \otimes \C [ t, t^{-1} ]$. Its bracket is induced by the relations $\left[ x \otimes u, y \otimes v \right] = \left[ x,y \right] \otimes uv$ for every elements $x,y \in \mathfrak{g}$ and $u,v \in \C [ t, t^{-1} ]$. Many natural infinite-dimensional representations of this algebra are projective representations (\textit{cf}. \cite{kac:bombaylectures}), signalling the relevance of a central extension of this Lie algebra. Recall that a \textit{central extension} of $\widetilde{\mathfrak{g}}$ by a Lie algebra $\mathfrak{a}$ is a Lie algebra $\mathfrak{q}$ containing $\mathfrak{a}$ in its center, \textit{i.e.} such that $\mf{a}$ has bracket $0$ with every element of $\mf{q}.$ We then have an exact sequence \be{exactSequenceCentralExt} 0 \longrightarrow \; \mathfrak{a} \longrightarrow \; \mathfrak{q} \longrightarrow \; \widetilde{\mathfrak{g}} \longrightarrow 0 \ee of Lie algebras. Using the notion of morphisms of central extensions (\textit{cf}. \cite[Eq. 1.3]{gar:loopGroups}), we can define the category of central extensions of $\mf{g}$. An initial object in this category is called a \textit{universal central extension}\footnote{If $\mf{q}$ denotes such an initial object, for any central extension $\mf{q}'$ of $\mf{g}$, there exists a unique morphism $\mf{q} \longrightarrow \mf{q}'$ of central extensions.}. If $\widetilde{\mf{g}}$ is perfect, \textit{i.e.} if we have $[\widetilde{\mf{g}}, \widetilde{\mf{g}}]=\widetilde{\mf{g}}$, or in other words if $\widetilde{\mf{g}}$ is generated by the brackets of its elements, such an initial object always exists up to unique isomorphism, using \cite[Prop. 1.7, Lem. 1.10]{gar:loopGroups}. It will be referred to as \emph{the} universal central extension of $\widetilde{\mf{g}}$.
            
    \begin{proposition}\cite[Thm. 2.36]{gar:loopGroups}
        The loop algebra $\widetilde{\mathfrak{g}}$ is perfect and its universal central extension $\gext$ is given by a one dimensional algebra $\mf{a} = \C \mathbf{c}$, with $\cc$ a central element. As a vector space, we have \be{univ:ext} \begin{array}{lll} \widehat{\mathfrak{g}} & = & \widetilde{\mathfrak{g}} \; \oplus \; \C \mathbf{c}, \end{array} \ee  where $\mathbf{c}$ is a central element, and the Lie bracket is  specified by linearly extending the following rule: \be{univ:ext:bracket} \begin{array}{lll} [ X \otimes t^n + \alpha \cc, \; Y \otimes t^m + \beta \cc ] & = & [X, Y] \otimes t^{n+m} \; + \; m \delta_{m, n} (X, Y) \,  \cc, \end{array} \ee for $X, Y \in \mf{g}$, $\alpha, \beta \in \C$ and where $(\cdot , \cdot)$ is the normalized Killing form introduced in \eqref{norm:kill}. Note that the inclusion of $\widetilde{\mathfrak{g}}$ into $\gext$ is not a morphism of Lie algebras, since it is not compatible with the brackets.
    \end{proposition}

\tpoint{Extensions by the derivation $\dd$} One does not obtain a satisfactory notion of roots for the Lie algebra~$\gext$ if one chooses to use  \be{h:hat}  \wh{\mf{h}} := \mf{h} \oplus \C \cc \ee as the analogue of a Cartan subalgebra. Another manifestation of this phenomenon is the fact that many of the interesting representations of $\gext$ actually possess a natural grading\footnote{The grading is usually an indication of another infinite dimensional Lie algebra called the Virasoro algebra, which also acts on the representation (\textit{cf}. \cite{segal:unitary}).}. For these reasons, it becomes natural to consider a further extension of $\gext$. Let $\dd: \gext \rr \gext$ be the derivation of $\gext$ which verifies \be{act:d} \begin{array}{llllllll} \dd|_{\mf{h} \oplus \C c} & = & 0 & \text{ and } & \dd( X \otimes t^n) & = & n X \otimes t^n & \text{ where } \; X \in \mf{g}, \; n \in \zee. \end{array} \ee  One can then consider the semi-direct sum (of Lie algebras, \cite[Ch. 1 \S 8]{bourbaki:1}) of $\gext$ with the one-dimensional Lie algebra $\C \dd$, and we denote this object as $\gext^e$. As a vector space, we have
		\be{ge:loop} 
		    \begin{array}{ccc}
                \gext^e & = & \gext \oplus \C \dd
            \end{array}, \ee and the Lie bracket is defined as  
        \be{ge:loop-bracket} 
            \begin{array}{ccc}
                \left[  x + \lambda \dd , y + \mu \dd \right] & = & \left[ x, y \right] + \lambda \dd \left( x \right) - \mu \dd \left( y \right)
            \end{array} . \ee Now the correct analogue of a Cartan subalgebra inside of $\gext^e$ is the $\ell+2$-dimensional vector space \be{cartan:e} \begin{array}{ccccc} \haff^e & := & \haff \oplus \C \dd & = & \C \cc \oplus \mf{h} \oplus \C \dd \end{array}. \ee 
            
\addtocounter{theorem}{1}

\tpoint{Completions for the $t$-adic topology} \label{subsubTadic}  In the definition of $\widehat{\mf{g}}$, we could have replaced the Laurent polynomials $\C[t, t^{-1}]$ with the formal Laurent series $\C((t))$. All of the constructions above go through without any change, and the resulting $t$-adically complete algebras will be denoted as $\gext_{\comp}$ and $\gext^e_{\comp} $. Note that changing the scalars to $\C((t))$ commutes with taking the semi-direct product with $\dd$, so the meaning of~$\gext^e_{\comp} $ is unambiguous.

\addtocounter{theorem}{1}

\tpoint{Affine Cartan matrices} \label{subsubAffineCartan}

Another useful point of view on the Lie algebra $\gext^e$ is its realization as an untwisted affine Kac--Moody algebra. To the $\ell \times \ell$ Cartan matrix $A$ introduced in \eqref{cartan-matrix}, one can attach the untwisted \emph{affine Cartan matrix}  \be{}\label{affineCartanMatrix} \begin{array}{llll}
            \displaystyle \wh{A}_{i j} & \displaystyle = & \displaystyle 2 \frac{\left( \alpha_i, \alpha_j \right)}{\left( \alpha_i, \alpha_i \right)} & \text{ for } \; i, j = 1, \ldots, \ell+1 
        \end{array}, \ee where we define $\alpha_{\ell+1} = -\alpha_0$ as the opposite of the highest root of $\mf{g}$ (\textit{cf}. \eqref{highest-root}). Note that we recover the Cartan matrix $A$ as the first $\ell$ rows and columns of $\wh{A}$. The matrix $\wh{A}$ is also an example of a \textit{generalized Cartan matrix} (\textit{cf}. \cite[Sec. 1.1]{kac:book}), and such matrices are classified into three main types: finite, affine, and indefinite, as explained in \cite[Thm. 4.3]{kac:book}. The affine ones have a further division into untwisted and twisted types, presented in \cite[Sec. 7.2, 8.3]{kac:book}, and $\wh{A}$ is actually an example of an \emph{untwisted affine type} matrix. It is also an example of a \emph{symmetrizable} generalized Cartan matrix (see \cite[Sec. 2.1]{kac:book}), in that there exists a diagonal matrix $D$ and a symmetric matrix $B$ such that we have $\wh{A} = D B$. The entries of $D$ are  $\epsilon_1, \ldots, \epsilon_{\ell+1}$, and given by $\varepsilon_i = 2/\left( \alpha_i, \alpha_i \right)$ here.
  
\addtocounter{theorem}{1}      

\renewcommand{\ad}{\mathrm{ad\,}}

\tpoint{Kac--Moody presentation} To the \emph{symmetrizable} generalized Cartan matrix $\wh{A}$, we can attach \footnote{Our presentation, in the symmetrizable case, is simplified by using the work of Gabber and Kac \cite{gabber-kac}, as explained in \cite[Thm. 9.11 and after]{kac:book}.} the corresponding \emph{Kac--Moody algebra $\mf{g}(\wh{A})$} as the quotient of the free Lie algebra generated by $3 \ell + 3$ symbols $e_1, \dots, e_{\ell+1}$, $f_1, \dots, f_{\ell+1}$, and $h_1, \dots, h_{\ell+1}$, together with the following relations:
    \[ \begin{array}{ccccccccccccccccccc}
            1. & \left[ h_i, h_j \right] & = & 0 & \quad & \hspace{-0.3em} 2. & \left[ e_i, f_j \right] & = & \delta_{i j} h_i & \quad & \hspace{-0.3em} 3. & \left[ h_i, e_j \right] & = & \wh{A}_{i j} e_j & \quad & \hspace{-0.3em} 4. & \left[ h_i, f_j \right] & = & - \wh{A}_{i j} f_j
        \end{array} \]
    \noindent for any integers $1 \leqslant i,j \leqslant \ell+1$, and
    \[ \begin{array}{ccccccccc}
            5. & \left( \ad e_i \right)^{-\wh{A}_{i j}+1} \left( e_j \right) & = & 0 & \qquad \qquad & 6. & \left( \ad f_i \right)^{-\wh{A}_{i j}+1} \left( f_j \right) & = & 0
        \end{array} \]
                
    \noindent for any integers $1 \leqslant i,j \leqslant \ell+1$ with $i \neq j$. Here, we have denoted by $\ad$ the adjoint representation, defined by $\left( \ad x \right) \left( y \right) = \left[ x,y \right]$ for $x, y \in \mathfrak{g} ( \wh{A} )$. The Kac--Moody algebra $\mf{g}(\wh{A})$ has a graded structure, which we will now describe. For nonnegative integers $n_1, \dots, n_{\ell+1}$ not all of which are zero, we set               
    \be{spanPositive}\begin{array}{ccc}
            \mathfrak{g} \left( n_1, \dots, n_{\ell+1} \right) & = & \spanVector_{\C} \left \{ \left[ e_{i_1}, \left[ e_{i_2}, \dots, \left[ e_{i_{r-1}}, e_{i_r} \right] \dots \right] \right] \right \}
        \end{array} ,\ee
    \noindent where each $e_i$ appears exactly $n_i$ times, and we set 
    \be{spanNegative}\begin{array}{ccc}
            \mathfrak{g} \left( -n_1, \dots, -n_{\ell+1} \right) & = & \spanVector_{\C} \left \{ \left[ f_{i_1}, \left[ f_{i_2}, \dots, \left[ f_{i_{r-1}}, f_{i_r} \right] \dots \right] \right] \right \}
        \end{array}\ee 
    \noindent where each $f_i$ appears $-n_i$ times. We further set
    \be{affineCartanSubalg}\begin{array}{ccccc}
            \mathfrak{h} ( \wh{A} ) & = & \mathfrak{g}_1 \left( 0, \dots, 0 \right) & = & {\C} h_1 \oplus \dots \oplus {\C} h_{\ell +1}
        \end{array}\ee
    \noindent To simplify the notation, we finally set $\mathfrak{g} \left( n_1, \dots, n_{\ell+1} \right) = 0$ for any other collection of integers $n_1, \dots, n_{\ell+1}$. Then from \cite[Sec. 3]{gar:loopGroups}  the Lie algebra $\mathfrak{g} ( \wh{A} )$ has a $\mathbb{Z}^{\ell+1}$-gradation, given by
    \be{ZlGrading}\begin{array}{ccc}
            \mathfrak{g} ( \wh{A} ) & = & \bigoplus\limits_{\left( n_1, \dots, n_{\ell+1} \right) \in \mathbb{Z}^{\ell+1}} \; \mathfrak{g}_1 \left( n_1, \dots, n_{\ell+1} \right).
        \end{array}\ee
Using the Chevalley basis of $\mf{g}$ reviewed in \S \ref{subsubChevalleyBasisIntegralForm}, we can now state the precise connection between $\mathfrak{g} ( \wh{A} )$ and the Lie algebra $\gext$ introduced in \S \ref{subsubLoop}.
               
    \begin{theorem}\cite[Thm. 3.14]{gar:loopGroups}
    \label{ThmIsoLoopKacMoody}
        There is an isomorphism of Lie algebras $\Psi: \gext \stackrel{\cong}{\longrightarrow} \mathfrak{g} ( \wh{A} )$, defined by \be{map:class} \begin{array}{lll} \left \{ \begin{array}{lll} \Psi \left( E_i \otimes t^0 \right) & = & e_i \\ \Psi \left( F_i \otimes t^0 \right) & = & f_i \\  \Psi \left( \av_i \right) & = & h_i \end{array} \right. & \text{  for } i=1, \ldots, \ell, \; \text{ and } & \left \{ \begin{array}{lll} \Psi \left( E_{- \alpha_0} \otimes t^1 \right) & = & e_{\ell+1} \\ \Psi \left( E_{ \alpha_0} \otimes t^{-1} \right) & = & f_{\ell+1} \\ \Psi \left( -\alpha_0^{\vee} + \cc \right) & = & h_{\ell+1} \end{array} \right. \end{array} . \ee
    \end{theorem}
            
\tpoint{Derivation in the Kac--Moody presentation} \label{ss:deg-der-KM} Recall that we introduced earlier the semi-direct sum of $\widehat{\mf{g}}$ with the operator $\dd$, which we denoted by $\gext^e$. We have a $\mathbb{Z}$-gradation of $\mathfrak{g} ( \wh{A} )$
            
            \begin{equation}
                \label{GradationKacMoodyZ}
                \begin{array}[t]{ccc}
                    \mathfrak{g} ( \wh{A} )_j & = & \bigoplus\limits_{\left( n_1, \dots, n_{\ell} \right) \in \mathbb{Z}^{\ell}} \; \mathfrak{g}_1 \left( n_1, \dots, n_{\ell}, j \right), 
                \end{array}
            \end{equation}
            
            \noindent  and define a linear map 
    \be{D:km} \begin{array}{ccccc}
        D & : & \mathfrak{g} ( \wh{A} )_j & \longrightarrow & \mathfrak{g} ( \wh{A} )_j \\[0.5em]             
         && x & \longmapsto & j x
    \end{array} .  \ee
	One may check that under the isomorphism $\Psi$ of Theorem \ref{ThmIsoLoopKacMoody}, the action of the map $\dd$ is sent to that of $D$ and it is also a derivation of $\mathfrak{g} ( \wh{A} )$. Let us now define the semi-direct product \be{gaff:ext} \begin{array}{lllll} \gaff & := &  \mathfrak{g} ( \wh{A} ) \rtimes \C D & = & \mathfrak{g} ( \wh{A} ) \oplus \C D. \end{array} \ee Our comments above indicate that the isomorphism of Theorem \ref{ThmIsoLoopKacMoody} extends to the extensions by $\dd$ and $D$ respectively. We continue to denote this isomorphism by the same name \be{ext:KM} \begin{array}{lllll} \Psi & : & \mathfrak{g}^e ( \wh{A} )  & \stackrel{\cong}{\longrightarrow} & \gext^e. \end{array} \ee  We write $\mf{h}^e(\wh{A}):= \mf{h}(\wh{A}) \oplus k D$ and note that it is isomorphic to the Cartan subalgebra $\wh{\mf{h}}^e$ introduced in \eqref{cartan:e} above. One gets the analogue for $\mathfrak{g}^e ( \wh{A} )$ of the $\mathbb{Z}$-gradation induced by \eqref{GradationKacMoodyZ} by setting \be{gradationKacMoodyCartan} \begin{array}{ccccc}
                        \mathfrak{g} ( \wh{A} )_0 & = & \mathfrak{h}^e ( \wh{A} ) & = & \mathfrak{h} ( \wh{A} ) \oplus k \mathbf{d}
                    \end{array}\ee
    \noindent without changing the other spaces in the gradation. 
    
    \vspace{0.15in}

\addtocounter{theorem}{1}

\subsubsection{Completion of the Kac--Moody algebra}
\label{subTwoKacMoodyAlgebras}
        
In \S \ref{subsubTadic}, we completed the Lie algebra $\gext$ to $\gext_{\comp}$. To reproduce this procedure in the Kac--Moody context, recall the $\zee$-gradation introduced in (\ref{GradationKacMoodyZ}) introduced above. With respect to it, we define the \textit{norm} of  $x \in \mathfrak{g} ( \wh{A} )$ to be
  \be{def:norm} \begin{array}{ccc}
                        \Big\Vert \sideset{}{}\sum\limits_{j \in \mathbb{Z}} x_j \Big\Vert & = & e^{-j_0}
                    \end{array} \ee
                
                \noindent where only finitely many components $x_j$ of $x$ are non zero, and we have $j_0 \,= \, \min \left \{ j \in \mathbb{Z}, \; x_j \neq 0 \right \}$. The completion of $\mathfrak{g} ( \wh{A} )$ is denoted by $\mathfrak{g} ( \wh{A} )_\comp$.
            From \cite[Rmk. 2.8]{gar:loopGroups},  the isomorphism from Theorem \ref{ThmIsoLoopKacMoody} extends to an isomorphism of Lie algebras
            \be{ext:isom-comp} \begin{array}{ccccc}
                        \mathfrak{g} ( \wh{A} )_{\comp} & \simeq & \mathfrak{g} \otimes k (( t )) \oplus k \mathbf{c} & =: & \widehat{\mathfrak{g}}_\comp
                    \end{array} . \ee
                    \noindent If one considers in \eqref{ext:isom-comp} the gradation introduced in \eqref{gradationKacMoodyCartan}, the natural extension of the norm from \eqref{def:norm} to $\mf{g}^e(\wh{A})$ allows us to define a completion denoted as $( \mathfrak{g}^e ( \wh{A} ) )_{\comp}$. One can then verify that we have \be{commute:comp-der} \begin{array}{lll} ( \mathfrak{g}^e ( \wh{A} ) )_{\comp} & = & ( \mathfrak{g} ( \wh{A} ) )_{\comp} \; \rtimes \; \C D \end{array} \ee where the action of $D$ on $( \mathfrak{g} ( \wh{A} ) )_{\comp}$ considered in \eqref{D:km} extends to a derivation of the completion. Note that the above is also isomorphic to $\gext^e_{\comp}$ introduced in \S \ref{subsubTadic}.

\tpoint{Notational convention} \label{notationLoopAlgebras} From now on, we refer to $\gext^e$ as either the loop or affine Kac--Moody algebra, and shall no longer write $\gaff$. Similarly, we write $\gext^e_{\comp}$ even when we might mean $\gaff_{\comp}$. In light of the explicit isomorphisms given above, we hope this will not lead to any confusion.

 \tpoint{Roots of the extended loop algebra $\gext^e$} Let us now describe the root system\footnote{We use the term `root system' informally, since the standard axioms for root systems and their references, \textit{e.g.}  \cite[Ch. VI]{bourbaki:456}, pre-date Kac--Moody theory by a few years and so do not actually admit the root systems of affine Lie algebras as examples. Nonetheless, many of the properties discussed in \emph{op. cit.} do carry over.} of $\gext^e$. Recall that we have described the root system $\rts \subset \mf{h}^*$  and its basis $\Pi = \{ \alpha_1, \ldots, \alpha_{\ell} \}$ of simple roots  in \S \ref{subsubFinRts}.  For each linear form $\lambda \in (\haff^e)^*$ we define the space 
                \be{root:sp:aff} \begin{array}{lll} (\gext^e)^{\lambda} & = & \{ X \in \gext^e \mid [H, X] \; = \; \lambda (H) X \; \text{ for all } \; H \in \haff^e \}, \end{array} \ee and refer to those elements $\lambda \in (\haff^e)^*$ such that we have $(\gext^e)^{\lambda}\neq 0$ as the \emph{roots} of $\gext^e.$ Denote the set of roots as $\rtsaff$, and, as before, for any root $a \in \rtsaff$, we refer to $(\gext^e)^a$ as the corresponding \textit{root space}. Much of the theory of the finite root system $\rts$ carries over to $\rtsaff$, but a new feature appears, namely one has `imaginary' roots such that the dimension of the corresponding root space is no longer $1$.\footnote{In the case of affine Kac--Moody algebras, the multiplicites of these imaginary roots is always $\ell$, but computing such multiplicities in a general Kac--Moody algebra remains quite mysterious.}
                
                 Let us now explicitly describe the roots $\rtsaff$. First, define the (minimal, positive) \textit{imaginary root} $\iota \in ( \haff^e )^{*}$ using the decomposition  \eqref{cartan:e} as follows:                 
                \be{def:iota} \begin{array}{ccc}
                        \iota \left( x + \lambda \mathbf{d} \right) & = & \lambda \text{ for } x \in \C \cc \oplus \mf{h}.
                    \end{array}  \ee 
                     For any integer $1 \leqslant i \leqslant \ell$, we extend $\alpha_i \in \Pi$  to an element $a_i \in (\haff^{e})^*$  defined as  \be{ai:class} \begin{array}{llllllllll} a_i |_{\mf{h}} & = & \alpha_i, \; \; \; & a_i(\cc) & = &0, \; \; \; & \text{ and } & a_i(\dd) & = & 0. \end{array} \ee As every element in $\rts$ is a linear combination of $\alpha_1$, \ldots, $\alpha_{\ell}$, we use the construction above to extend any root $\alpha \in \rts $ to an element of $(\haff^e)^*$. Unless explicitly mentioned (\textit{e.g.} in the case of $a_i$ for $i=1, \ldots, \ell$), we shall denote this extension by zero of an element of $\rts$ into an element of $(\haff^e)^*$ by the same name. Finally, let us introduce the element $a_{\ell+1} \in (\haff^{e})^*$ defined by 
                \be{a:new} \begin{array}{lll} a_{\ell+1} & = &  - \alpha_0 + \iota . \end{array} \ee Note that the elements $a_1, \dots, a_{\ell+1} \in (\haff^e)^*$ are now linearly independent, a fact which would not be true if we regarded these elements as functionals on $\haff$. 
                                
               With these preliminaries, we can now describe the set $\rtsaff.$ It has a decomposition into positive and negative roots $\rtsaff = \rtsaff_+ \sqcup \rtsaff_-$ where we have $\rtsaff_- = - \rtsaff_+$ and     
                \be{rts:desc}  
                        \begin{array}{lll} \rtsaff_+ & = &  \{ \alpha + n \iota \mid \alpha \in \rts_+, n \in \zee_{\geq 0} \} \; \sqcup \; \{ -\alpha + n \iota \mid \alpha \in \rts_+, n \in \zee_{\geq 1} \} \; \sqcup \; \{ n \iota \mid n \in \zee_{\geq 1}  \} . \end{array} \ee The set of \emph{imaginary roots} $\rtsaffim$ and \emph{Weyl (or real) roots} $\rtsaffre$ are then defined as \be{real:im} \begin{array}{lllllll} \rtsaffim & := & \{ n \iota \mid n \in \zee_{\neq 0} \} & \; \text{ and } \; & \rtsaffre & := & \rtsaff \setminus \rtsaffim \end{array} \ee  respectively. We also define the positive/negative real roots as $(\rtsaffre)_{\pm}:= \rtsaffre \cap \rtsaff_{\pm}$ and similarly define the positive and negative imaginary roots $(\rtsaffim)_{\pm}$. The set $\affpi := \{ a_1, \ldots, a_{\ell+1} \} $ plays the role of a basis for $\rtsaff$ in that every $a \in \rtsaff$ can be written as a $\zee$-linear combination of elements from $\affpi$ with either all non-negative or all non-positive coefficients. For future reference, we also define the \emph{root lattice} $\widehat{\mathscr{Q}}$ as \be{rt:lattice} \begin{array}{lll} \rtl & := & \zee a_1 + \zee a_2 +  \cdots + \zee a_{\ell+1} . \end{array} \ee We make the natural definition for $\widehat{\mathscr{Q}}_+$, as $\mathbb{N}$-linear combinations.

\tpoint{Invariant bilinear form}  \label{subsubKacForm} Recall from \S \ref{subsubAffineCartan} that we have defined the entries $(\epsilon_1, \ldots, \epsilon_{\ell+1})$ of the diagonal matrix $D$ which symmetrizes $\wh{A}$, \textit{i.e.} such that we have $\wh{A} = D B$ with $B$ symmetric. Using the entries of this matrix, we can now define a symmetric bilinear form $(\cdot, \cdot)$ on $\haff^e$ by the conditions \be{kac-form} \begin{array}{lllllll} (h_i, h) & = & a_i(h) \epsilon_i & \text{ for } i \in \lb 1, \ell+1 \rb \text{ and } h \in \haff^e \text{, and by } & (\dd, \dd) & = & 0. \end{array} \ee Note that the first condition implies that we have \be{norm:D} \begin{array}{lllllll} (h_i, \dd) & = & 0 & \text{ for } i \in \lb 1, \ell \rb \text{, \; and } & (h_{\ell+1}, \dd) & = & 1. \end{array} \ee This form is non-degenerate (\cite[Lemma 2.1b]{kac:book}) and hence induces an isomorphism $\nu: \haff^e \rr (\haff^e)^*$. Under this isomorphism, we have $\nu(h_i) = \epsilon_i a_i$ and the corresponding form on $(\haff^e)^*$, again denoted by $(\cdot, \cdot)$, satisfies \be{dual-kac-form}  \begin{array}{llllllll} (a_i, a_j) & = & \wh{A}_{ij} / \epsilon_i & = & b_{ij} & \text{ for } \; i, j \in \lb 1, \ell+1 \rb. \end{array} \ee One can also check that we have \be{iota:norm} \begin{array}{llllllll} (\iota, \iota) & = & 0 & \text{ and } & (\iota, a_i) & = & 0 & \text{ for } i =1, \ldots, \ell. \end{array} \ee This form agrees with the normalized Killing form on $\mf{h}^*$ constructed in \S \ref{subsubKillingCoroots}.

        \newcommand{\affW}{\wh{W}}
        \newcommand{\haffd}{(\haff^e)^*}

\tpoint{Coroots and integral weights} \label{subsubCoroots}  Let us now work with the Kac--Moody presentation to define the coroots. Recall that we have introduced, within this presentation, elements $h_1, \ldots, h_{\ell+1}$ of $\mf{h}(\wh{A})$. We define \be{corts:classical} \begin{array}{llll} \av_i & = & h_i & \text{ for } \; i=1, \ldots, \ell+1. \end{array} \ee Some care has to be taken to introduce the coroots attached to a general element $\varphi \in \rtsaff$. We proceed first by setting $h'_{a_i}:= \epsilon_i^{-1} h_i$, and then by extending this into $h'_{\varphi}:= \sum_{i=1}^{\ell+1} c_i h'_{a_i}$ for any $\varphi \in \rtsaff$ written as $\varphi = \sum_{i=1}^{\ell+1} c_i a_i$. Finally, we can define the coroot $\varphi^{\vee}$ attached to $\varphi$ as \be{def:crt} \begin{array}{lll} \varphi^{\vee} & = & \displaystyle \frac{2}{(\varphi, \varphi)} \, h'_{\varphi}. \end{array} \ee We call $\rtsaff^{\vee}:= \{ \varphi^{\vee} \mid \varphi \in \rtsaff \}$ the set of affine coroots (it is the set of roots of the affine, possibly \emph{twisted}, Kac--Moody algebra with Cartan matrix $\leftidx^T A$.) The set $\wh{\Pi}^{\vee}:= \{ \av_1, \ldots, \av_{\ell+1} \}$ is then a `basis' of $\wh{\mathscr{R}}^{\vee}$, in the same sense that $\Pi$ was a basis of $\rtsaff$. Using the coroots, we can define the lattice of \emph{integral weights} as \be{lat:wts} \begin{array}{lll} \Lambda & := & \{ \lambda \in (\haff^e)^* \mid \la \lambda, \av_i \ra \in \zee \text{ for } i \in \lb 1, \ell+1 \rb \text{ and } \la \lambda , \dd \ra \in \zee \} . \end{array} \ee Note that we have $\wh{\mathscr{Q} } \subset \Lambda$, and we introduce the \emph{dominance order} $\preccurlyeq$ on $\Lambda$ as follows: for any $\lambda, \mu \in \Lambda$, we write $\lambda \preccurlyeq \mu$ if we have $\mu - \lambda \in \wh{\mathscr{Q}}_+.$ We also say that an element $\lambda \in \Lambda$ is \emph{dominant} if we have $\la \lambda, \av_i \ra \geq 0$ for all $i \in \lb 1, \ell+1 \rb,$ and the set of such elements is written as $\Lambda_+.$  Of particular importance to us will be the \emph{fundamental weight} $\Lambda_{\ell+1} \in \Lambda$, defined by the conditions \be{fund:wt} \begin{array}{lllllll} \la \Lambda_{\ell+1} , \av_i \ra & = & \iota_{i,\ell+1} & \text{for } \; i \in \lb 1, \ell+1 \rb \; \text{ and } & \la \Lambda_{\ell+1}, \dd \ra  & = & 0. \end{array} \ee  One then notes that \be{hdual:dec} \haffd  \; \text{ has basis } \; a_1, a_2, \ldots, a_{\ell+1}, \Lambda_{\ell+1}. \ee As we can write $a_{\ell+1}$ in terms of the classical roots $a_1, \ldots, a_{\ell}$ and $\iota$ we also have \be{hdual:2} \begin{array}{lll} \haffd & = & \C \iota \oplus \mf{h}^* \oplus \C \Lambda_{\ell+1}. \end{array} \ee  With respect to the inner product introduced in \S \ref{subsubKacForm}, one may check that we have \be{Lambda:norm} \begin{array}{ccc} (\Lambda_{\ell+1}, \Lambda_{\ell+1}) = 0, & (a_i, \Lambda_{\ell+1})=0 \; \text{ for } \; i \in \lb 1, \ell \rb, & (a_{\ell+1}, \Lambda_{\ell+1}) =1. \end{array} \ee 
For $\mu \in \haffd$, we often write $\overline{\mu}$ for its projection onto $\mf{h}^*$ and record a useful formula (\textit{cf}. \cite[(6.2.6)]{kac:book}) \be{lambda:proj} \begin{array}{lll} \lambda & = & \la \lambda, \cc \ra \Lambda_{\ell+1} +  \overline{\lambda} + (\lambda , \Lambda_{\ell+1}) \iota. \end{array} \ee

\tpoint{The affine Weyl group } 
\label{subsubAffineWeylGroup}
    
For $i \in \llbracket 1, \ell + 1 \rrbracket$, define the reflection            
    \be{defReflexions}
        \begin{array}{lllll} r_i  & : &  (\haff^e)^*  & \longrightarrow  & (\haff^e)^* \\[0.4em]  && \lambda & \longmapsto & \lambda - \lambda \left( \av_i \right) a_i \end{array}
       .\ee The group generated by these reflexions is denoted by $\wh{W} $ and is called the \textit{affine Weyl group}. It also has a presentation as $\wh{W} = W \rtimes \mathscr{Q}^{\vee}$ where $\mathscr{Q}^{\vee}$ is the coroot lattice (\textit{cf.} \cite[\S 6.5]{kac:book}), though we will not need this explicitly here. The group $\wh{W}$ is also a Coxeter group, with Coxeter generators $\wh{S}:= \{ r_i \mid i \in \lb 1, \ell+1 \rb \}$. Let us note that the form $(\cdot, \cdot)$ on $(\haff^e)^*$ introduced in the previous section is $\wh{W}$-invariant (\textit{cf}. \cite[Prop. 3.9]{kac:book}) and that every element of $\wh{W}$ fixes $\iota$, whereas every element in $\rtsaffre$ is of the form $w a_i$ for some elements $w \in \affW$ and $a_i \in \wh{\Pi}$. For this reason, the elements in $\rtsaffre$ are sometimes referred to as the \textit{Weyl roots}.

\subsubsection{The Tits cone } 
\label{subsubTisCone}     
            
Although not explicitly needed, let us mention a crucial difference between the actions of $\wh{W}$ on $\haffd$ and of $W$ on $\mf{h}^*$: whereas every element $\lambda \in \mf{h}^*$ is $W$-equivalent to an element in the dominant cone \be{dominantCone} \begin{array}{lll} \mc{C} & := & \{ \lambda \in \mf{h}^* \mid \la \lambda, \alpha^{\vee}_i \ra \geq 0 \text{ for } i \in \lb 1, \ell \rb \}, \end{array} \ee this is no longer true in $\haffd$. Instead, we need to introduce the \textit{Tits cone} \be{tits:cone} \begin{array}{lllllll} X & := & \bigcup_{w \in \wh{W}} \, w (\wh{\mc{C}}), & \text{ where } & \wh{\mc{C}} & := & \la \lambda \in \haffd \mid \la \lambda, \av_i \ra \geq 0 \text{ for } i \in \lb 1, \ell+1 \rb \}. \end{array} \ee  Note that the projection of $\wh{\mc{C}}$ onto $\mf{h}^*$, following the decomposition \eqref{hdual:2}, is compact.  As it turns out, the Tits cone can be described explicitly as \be{tits:explicit} \begin{array}{lll} X & = & \{ \lambda \in \haffd \mid \la \lambda, \cc \ra > 0 \}. \end{array} \ee  We also remark here that at the group level, which will be presented later, working with the condition $|\tau|<1$ or $|\tau|>1$ is tantamount to working inside or outside the Tits cone.

\newcommand{\affb}{\wh{\mf{b}}}
\newcommand{\affu}{\wh{\mf{u}}}
\newcommand{\ueg}{\mc{U}}
\newcommand{\mult}{\mathrm{mult}}

\subsection{Highest weight representations}
\label{subGarlandRep} 
Let us review some basic features of the theory of highest weight representations of $\gext^e$ over $\C$ here. Unless otherwise mentioned, our base field will be $\C$ throughout this section, and later, after introducing integral forms, we extend this construction to an arbitrary field.

\tpoint{Weight modules} \label{rep:basics} A representation $V$ of $\gext^e$ is called is called a \emph{weight module} if \begin{enumerate} 
                \item we have  $V = \bigoplus_{\mu \in \haffd} \, V_{\mu}$, where we set \be{wt:spaces} \begin{array}{lll} V_{\mu} & = & \{ v \in V \mid H.v= \lambda(H) v \text{ for all } H \in \haff^e\}; \end{array} \ee 
                \item each $V_{\mu}$ is finite-dimensional. \end{enumerate} 
The set of all elements $0 \neq \mu \in \haffd$ such that we have $V_{\mu} \neq 0$ is denoted by $\wts(V)$ and its elements are called the \textit{weights} of $V$. We define the \emph{multiplicty of the weight} $\mu \in \wts(V)$ as \be{}\label{mult:wt} \begin{array}{lll} \mult(\mu) & := & \dim_{\C} V_{\mu}. \end{array} \ee  For a weight module $V$ we define its \textit{formal character} as the sum \be{}\label{formal:char} \begin{array}{lll} \chi_V & := & \sideset{}{}\sum\limits_{\mu \in \wts(V)} \mult(\mu) \, e^{\mu} \end{array} \ee where $e^{\mu}$ is a formal symbol which lives in some algebra (\textit{cf}. \cite[\S 9.7]{kac:book} for the precise details), satisfying the condition $e^{\mu} e^{\nu} = e^{\mu + \nu}$. Let us recall that a representation $V$ is called \emph{integrable} if the elements $e_i, f_i \in \gext^e$ act \emph{locally nilpotently} for each $i \in \lb 1, \ell+1 \rb$, \textit{i.e.} if \be{}\label{loc:nilp} \text{ for any } v \in V \text{ there exists an integer } n \text{ such that we have } e_i^n v = f_i^n v =0. \ee To exponentiate certain operators and define a group, we need this condition.

\addtocounter{theorem}{1}

\tpoint{Enveloping algebra} \label{subsubEnvelopingAlgebras} For any complex Lie algebra $\mf{q}$, its \textit{universal enveloping algebra}  $\mathcal{U}(\mathfrak{g})$ is defined as the quotient of the tensor algebra \be{tensor:alg} \begin{array}{lll} T(\mf{q}) & := & \sideset{}{}\bigoplus\limits_{n \geq 0}  \; \; \underbrace{(\mf{q} \otimes \cdots \otimes \mf{q})}_{n-\text{times}} \end{array} \ee by the smallest ideal containing the quadratic relation $x \otimes y - y \otimes x = [x, y] $ for $x, y \in \mf{q}$. Note that if $\mf{q}_1 \subset \mf{q}$ is a subalgebra, one has a morphism of algebras $\mc{U}(\mf{q}_1) \rr \mc{U}(\mf{q})$. The importance of these algebras stems from the fact that a representation of $\mf{q}$ is the same as an $\mc{U}(\mf{q})$-module. 

\addtocounter{theorem}{1}

\tpoint{Unipotent and Borel subalgebras} \label{subsubUnipotentSubalgebras} Define the following important subalgebras of $\gext^e$,                
    \be{unip} \begin{array}{ccccc}  
        \wh{\mathfrak{u}}_{\pm}   =  \bigoplus\limits_{\alpha \in \rtsaff_{\pm} } \; ( \gext^e )^{\alpha},  & \wh{\mf{b}}_{\pm}  =  \wh{\mathfrak{h}}  \oplus \wh{\mathfrak{u}}_+, & \text{ and } & \wh{\mf{b}}_{\pm}^{e}:= \wh{\mf{b}} \oplus \C \dd.  \end{array}
        \ee Usually we omit the subscript `+' if context makes it clear we are speaking of the positive case. Note that we have a triangular decomposition \be{triangle} \begin{array}{lll} \gext^e & = & \affu_- \oplus \haff^e \oplus \affu_+, \end{array} \ee which corresponds by the Poincare--Birkhoff--Witt theorem to a decomposition \be{PBW} \begin{array}{lll} \ueg(\gext^e) & = & \ueg(\affu_-) \; \otimes_{\C} \; \ueg(\haff^e) \; \otimes_{\C} \; \ueg(\affu_+). \end{array} \ee  

\addtocounter{theorem}{1}

\tpoint{Irreducible highest weight modules $V^{\lambda}$}  \label{subsubIrreducible} 
To each element $\lambda \in \haffd$ we attach the \textit{Verma module} \be{verma} \begin{array}{lll} M \left( \lambda \right) & = & \mc{U}(\gext^e) \otimes_{\mc{U}(\affb^e)} \C_{\lambda}, \end{array} \ee where $\C_{\lambda}$ is the one-dimensional $\affb^e$-module on which $\wh{\mf{h}}^e$ acts by $\lambda$, \textit{i.e.} such that we have $H.v = \lambda(H) v$ for every $H \in \haff^e$ and $v \in \C_{\lambda},$ and on which $\affu_+$ acts trivially. Note that, by (\ref{PBW}), we have \be{M: pbw} \begin{array}{lll} M(\lambda) & \cong & \ueg(\affu_-) \otimes_{\C} \C_{\lambda}, \end{array} \ee and from this one verifies that $M(\lambda)$ is a weight module.  Any non-zero quotient of $M(\lambda)$ is called a \textit{highest weight module} with highest weight $\lambda$ since every weight $\mu \in \wts(M(\lambda))$ satisfies $\mu \preccurlyeq \lambda$, where we recall that $\preccurlyeq$ refers to the dominance order introduced in \S \ref{subsubCoroots}. One knows (\textit{cf}. \cite[Lem. 2.1.2]{kumar}) that $M(\lambda)$ has a unique proper maximal submodule $R_{\lambda}$ such that the quotient $V^{\lambda}:= M(\lambda) / R_{\lambda}$ is irreducible as a representation. The weight space of weight $\lambda$ has dimension~$1$, in both $M(\lambda)$ and any of its quotients. We pick a non-zero vector $v_{\lambda}$ in this weight space and call it a \emph{highest weight vector}. We restrict our attention to the modules $V^{\lambda}$ and one actually knows (\cite[Cor. 2.1.8]{kumar}) that the map $\lambda \mapsto V^{\lambda}$ is a bijective correspondence between the set of dominant, integral weights $\Lambda_+$ (\textit{cf}. right after \eqref{lat:wts}) and the isomorphism classes of integrable (\textit{cf.} \S \ref{rep:basics}), irreducible, highest weight modules of $\gext^e$. For such a representation, we denotes its set of weights as $\wts_{\lambda}$ for short.

\addtocounter{theorem}{1}

\tpoint{Depth and level}  \label{subsubDepthLevel} \label{conv:hw} Fix $\lambda \in \Lambda_+$ such that we have $\la \lambda, \cc \ra > 0.$ Such weights are said to be of \emph{positive level}, and we write $V:= V^{\lambda}$ for the corresponding irreducible, highest weight representation with weight lattice $\wts_{\lambda}$. In the remainder of this subsection, we collect a number of facts about this set $\wts_{\lambda}$ which we need for our future work. To state these, we introduce the \emph{depth} $\depth(\mu)$ and \emph{level} $\lev(\mu)$ of a weight $\mu \in \wts_{\lambda}$, by: since $\mu \preccurlyeq \lambda,$ we have $\mu   =  \lambda - m_1 a_1 - \cdots m_{\ell+1} a_{\ell+1}$ for some non-negative integers $m_i$. We then set
    \be{depth:level}             
        \begin{array}{lllllll} \depth (\mu)  & = &  \sum\limits_{i=1}^{\ell+1} \; m_i  & \text{ and } &  \lev(\mu) & = & m_{\ell+1}. \end{array} \ee 

\addtocounter{theorem}{1}

\tpoint{Linear-quadratic inequality} We now state an important inequality relating the classical (\textit{i.e.} projection onto $\mf{h}^*$) and imaginary parts (\textit{i.e.} projection onto $\C\mf{\iota}$) of a weight from a highest weight representation. Fix $\lambda \in \Lambda_+$ of positive level $p = \la \lambda, \cc \ra.$ For any $\mu \in \wts_{\lambda}$, which we may write as in (\ref{lambda:proj}) \be{mu:form} \begin{array}{lll} \mu & = & p \Lambda_{\ell+1} + \overline{\mu} + n \iota. \end{array} \ee From \cite[Proposition 2.12d]{kp}, we know that there exists some $a_0 \in \{ 1, 2 \}$ such that we have \be{pol-est} \begin{array}{lll} \vert \overline{\mu} \vert^2 - 2 a_0^{-1} \, p \, n & \leq & |\lambda|^2, \end{array} \ee where the norm on $\mf{h}^*$ is the restriction of $(\cdot, \cdot).$ Hence, there exists $C:=C(\rts, \lambda) > 0$ such that we have \be{lin-quad} \begin{array}{lll} \vert \overline{\mu} \vert & \leqslant & C \sqrt{n}. \end{array} \ee We may expand $\overline{\mu}$ in terms of the basis $\Pi$ as $\overline{\mu} = q_1 a_1 + \cdots + q_{\ell} a_{\ell}$ where $q_j \in \zee$ (they are not always positive). The inequality above then tells us, if we compare the norm $(\cdot, \cdot)$ with the sup-norm on the finite-dimensional space $\mf{h}^*$, that we also have $q_i \leq C \sqrt{n}$ for a (possibly different) constant $C:= C(\rtsaff, \lambda)$.

\tpoint{Maximal weights and multiplicities of weight spaces}  A \emph{maximal weight} $\mu \in \wts_{\lambda}$ is defined as a weight such that we have $\lambda + \iota \notin \wts_{\lambda}$, and we denote by $\max(\lambda)$ the set of all maximal weights. For every weight $\mu \in \wts_{\lambda}$, there exists a unique maximal weight $\eta \in \max(\lambda)$ and a unique non-negative integer $n'$ such that we have $\mu = \eta - n' \iota$, with $n' \leq \lev \mu$, using \cite[Prop. 12.5e]{kac:book}. From this and \cite[Thm. B, (4.24)]{kp}, we have the following.
\addtocounter{theorem}{1}
\begin{theorem} \label{thm:KPbound} There exist constants $A:= A(\lambda, \rtsaff) > 0$ and $C':= C'(\lambda, \rtsaff) > 0$ such that we have \be{mult:bd} \begin{array}{lll} \mult(\mu) & \leqslant & C' e^{A \sqrt{n}}\, , \end{array} \ee for any $\mu \in \wts_{\lambda}$ of level $n$. \end{theorem}  \noindent In fact a stronger asymptotic property in $n$ holds, though we shall not need this here. This estimate is proven by relating the Weyl--Kac character formula to the theory of theta functions and we refer to \emph{op. cit.} for the interesting details which go into its proof. Note the similarlity with the asymptotics of partition functions.

\addtocounter{theorem}{1}

\tpoint{Representation of the completed algebra}
\label{subsubRepCompletedAlgebra}
The previous paragraphs dealt with the representation of the extended affine Lie algebra $\gext^e \cong \mathfrak{g}^e ( \wh{A} )$. Using \cite[Prop 6.6]{gar:loopGroups}, for any Cauchy sequence $( x_n )_n$ of elements of $\mathfrak{g}^e ( \wh{A} )$, the sequence $x_n . v$ is ultimately constant. The representations studied above can thus be extended into representations of $\gext^e_{\comp} \cong \mathfrak{g}^e ( \wh{A} )_{\comp}$.

\subsection{Integral and Hermitian Structures}
\label{subIntegralFormGarlandRep}
        
In order to discuss the groups attached to the Kac--Moody algebras introduced in \S \ref{subSecLoopAlgebraUniversalCovering}, and in order to eventually describe the pro-Hermitian lattices which are the main subject of this work, we need to discuss certain integral structures on the Lie algebra $\gext^e$, its enveloping algebra $\ueg(\gext^e)$, and on its highest weight representations $V^{\lambda}.$ These were introduced by Garland in \cite{gar:loopAlgebras}. In the finite-rank case, the results are essentially due to Chevalley and Kostant (see \cite{kostant}  or \cite[Chap. 1-2]{steinberg:book} for an exposition). The situation is markedly more complicated in the loop setting owing in large part to the imaginary root spaces and their connection to the theory of symmetric functions.

\tpoint{Chevalley basis for $\gext$} \label{subsubChevalleyBasisLoop}   We begin by extending the basis $\mathcal{B}_{\mathrm{Ch}}(\mf{g})$ given at the end of \S \ref{subsubChevalleyBasisIntegralForm} to one for $\wh{\mf{g}}^e$, following \cite[Sec. 4]{gar:loopGroups} and \cite[Sec. 4]{gar:loopAlgebras}. Recall that we introduced the simple coroots $\av_i=h_i$ for every integer $i \in \llbracket 1, \ell+1 \rrbracket$ in \S \ref{subsubCoroots}. For every Weyl root $a \in \rtsaffre$ written as $a=\alpha + n \iota$, with $\alpha \in \rts(A)$ and $n \in \zee$, we first define the element in the root space $\wh{\mf{g}}^{a}$
        \be{xia}\begin{array}{ccc}
                        \xi_a & = & E_{\alpha} \otimes t^n
                    \end{array} .\ee
        \noindent Then, for every imaginary root $n \iota$, and every integer $i \in \llbracket 1, \ell \rrbracket$ we set
        \be{xiin}\begin{array}{ccc}
                        \xi_i \left( n \right) & = & H_i \otimes t^n
                    \end{array} .\ee
                
        \noindent A Chevalley basis of $\mathfrak{g} $ is then defined as the set 
                \be{aff:chev} \begin{array}{ccccccc}
                        \mathcal{B}_{Ch} ( \mathfrak{g} ( \wh{A} ) ) & = & \left \{ \av_1, \dots, \av_{\ell+1} \right \} & \cup & \left \{ \xi_a, \; a \in \rtsaffre(\wh{A}) \right \} & \cup & \left \{ \xi_i \left( n \right), \; n \in \mathbb{Z}^{\ast}, \; i \in \llbracket 1, \ell \rrbracket \right \}
                    \end{array}.  \ee From \cite[Thm. 4.12]{gar:loopAlgebras} we find that the structure constants of $\gext$ with respect to this basis are integers, and we refer to \emph{loc. cit.} for more precise and other favourable properties of this basis.

\tpoint{Integral forms of $\ueg(\gext)$ }  
The integral form $\ueg_{\zee}(\gext)$  is defined as the $\mathbb{Z}$-algebra generated by  \be{div-diff} \begin{array}{llll} \xi_a^{(p)} & := & \displaystyle \frac{\xi_a^p}{p!} & \text{  where } p\in \zee_{\geq0}, \; a \in \rtsaffre, \; \xi_a \text{ as in } \eqref{xia}. \end{array} \ee 

 \label{subsubIntegralFormZBasis}
\noindent Although we shall not need the rest of the results of this section in this work, we include them to give some further insight into the lattices we use to construct our pro-Hermitian systems in the next section. What we shall do, following \cite{gar:loopAlgebras}, is to now describe an explicit $\zee$-basis (as a module, not algebra) for this algebra.  Towards this end, we first introduce elements corresponding to $\ueg_{\zee}(\haff)$, namely for any integers $k \in \llbracket 1, \ell+1 \rrbracket$ and $r \in \mathbb{N}$, we set
                \be{}\label{} \begin{array}{ccc}
                        \displaystyle \binom{h_k}{r} & = & \displaystyle \frac{1}{r!} \, h_k \, \left( h_k - 1 \right) \, \dots \, \left( h_k - r + 1 \right)
                    \end{array} . \ee 
Next, for each $a \in \rtsaffre$, we have already introduced the elements which correspond to the one-dimensional real root spaces $\gaff^a$, namely the  $\xi_a^{(p)}$ from \eqref{div-diff}.  Finally, we need to define elements corresponding to the imaginary root space $\gaff^{ n \iota}$. We cannot take the elements $\xi_j(n)^p/p!$ since these elements are actually not even contained in $\ueg_{\zee}(\gext)$, nor do they preserve (via the adjoint action) the $\zee$-span of the Chevalley basis introduced above. Instead, as Garland discovered, one needs to introduce the following elements. In the polynomial algebra $\C \left[ X_1, X_2, \dots, X_n, \dots \right]$, define the polynomials $\Lambda_p (X) = \Lambda_p (X_1, \ldots, X_n, \ldots )$ in infinitely many variables by the formal identity \be{} \label{Lambda:gen-fns} \begin{array}{lll} \sideset{}{}\sum\limits_{p \geq 0} \Lambda_p(X) \, z^p & = & \exp\left( \sideset{}{}\sum\limits_{j \geq 1} \frac{X_j}{j} \, z^j \right). \end{array} \ee Using these polynomials we can define the elements \be{} \begin{array}{llll} \Lambda_p(r, j) & := & \Lambda_p(\xi_j(r), \xi_j(2r), \ldots) & \text{ for } p \geq 0, \; r \in \zee_{\neq 0}, \; \, j \in \lb 1, \ell \rb. \end{array} \ee Note that we have $\Lambda_0=1$ in our convention, meaning there is a shift of index with respect to the convention of \cite[p.500]{gar:loopAlgebras}. One checks also that we have $\Lambda_p(j) \in \ueg_{\zee}(\gext)$.

\addtocounter{theorem}{1}

Fix a total order on $\mc{B}_{\mathrm{Ch}}(\gext)$. Consider a sequence of non-negative integers \be{ints:mon} \begin{array}{llll} (p_a)_{a \in \rtsaffre}, & q(r, j)_{r \in \zee \neq 0} \; \text{ for } \; j \in \lb 1, \ell \rb, & \text{ and } & (\lambda_k)_{k=1, \ldots, \ell+1} \end{array} \ee such that almost all of the $p_a$ and $q(r, j)$ are zero.  With respect to such a family, we define a  monomial by taking the products with respect to our fixed order on $\mc{B}_{\mathrm{Ch}}(\gext)$ of the elements
\be{mon:seq} \begin{array}{lll} \displaystyle \frac{\xi_a^{p_a}}{p_a!}, &  \displaystyle \Lambda_{q(r, a_j)}(r, j), & \displaystyle \binom{h_k}{\lambda_k}. \end{array} \ee 
  
\begin{proposition}\cite[Thm 5.8]{gar:loopAlgebras}
    \label{PropZBasisUAtilde}
        The monomials as above form a basis of $\ueg_{\zee}(\gext)$, seen as a $\mathbb{Z}$-module.
\end{proposition}

\tpoint{Integral forms for representations}\label{SubSubIntFormsRep} Let $\lambda \in \Lambda_+$ be a dominant weight and $V^{\lambda}$ be the corresponding highest weight representation of $\gext^e$ introduced in \S \ref{subsubIrreducible}. Let $v_{\lambda}$ be a highest weight vector of $V^{\lambda}$, and define     
     \be{}\label{Vzee:form} \begin{array}{ccc}
                        V_{\mathbb{Z}}^{\lambda} & = & \mathcal{U}_{\mathbb{Z}} (\gext  ) \cdot v_{\lambda}.
                    \end{array}  \ee

\addtocounter{theorem}{1}

\tpoint{Admissible bases} \label{subsubAdmissible}A basis $\mc{B}$ of $V^{\lambda}$ is said to be \emph{admissible} if we have \be{adm:basis} \begin{array}{lll} \mc{B} & = & \bigsqcup\limits_{\mu \in \wts_{\lambda} } \mc{B} \; \cap \; V_{\mu}^{\lambda}. \end{array} \ee Not only do such bases exist, but in fact from \cite[Theorem 11.3]{gar:loopAlgebras}, we may choose one whose $\zee$-span is stable under $\ueg_{\zee}(\gext)$, where the vectors corresponding to a weight space are placed consecutively, and the weight spaces are ordered by increasing depth (see \S \ref{subsubDepthLevel}). Fix such a basis (called an \emph{integral admissible coherently ordered basis}). From \cite[Prop 11.7]{gar:loopAlgebras}, it follows that for such a basis, we have \be{} \begin{array}{lllllllll} V^{\lambda}_{\zee} & = & \bigoplus\limits_{\mu \in \wts_{\lambda}} \, V_{\zee, \mu}^{\lambda}, & \text{ where } & V^{\lambda}_{\zee, \mu} & := & V^{\lambda}_{\zee} \cap V_{\mu}^{\lambda} & = & \zee[ \mc{B}_{\lambda} \cap V^{\lambda}_{\mu} ]. \end{array} \ee 

\addtocounter{theorem}{1}

 \tpoint{Representations $V_k^{\lambda}$}\label{SubSubRepVk} For any field $k$, we can now define $V_k^{\lambda}:= V_{\zee}^{\lambda} \otimes_{\zee} k$, which is a representation of $\gext_k := \gext_{\zee} \otimes_{\zee} k$ with $\gext_{\mathbb{Z}}$ denoting the $\zee$-span of the Chevalley basis $\mc{B}_{\mathrm{Ch}}(\gext)$. Also, note that the element~$\dd \in \gext^e$ acts on $v \in V^{\lambda}_{\zee, \mu}$ as scalar multiplication by $\la \mu, \dd \ra \in \zee$. Hence, the subalgebra \be{} \begin{array}{lll} \gext^e_{\zee} & = & \gext_{\zee} \; \oplus \; \zee \dd, \end{array} \ee of $\gext^e$ also acts on $V^{\lambda}_{\zee}$.
 
 \addtocounter{theorem}{1}

\tpoint{Hermitian structure on $V^{\lambda}$}
\label{subHermitianStructure}
\label{subsubHermitianStructure}

Fix a $\zee$-form $V^{\lambda}_{\zee}$ as in \eqref{Vzee:form}. Then we recall  (see \cite[Thm. 12.1]{gar:loopAlgebras}) that $V^{\lambda}_{\mathbb{C}} = V^{\lambda}_{\mathbb{Z}} \otimes_{\mathbb{Z}} \mathbb{C}$ admits a \emph{positive-definite} Hermitian inner product $\{ \cdot , \cdot \}$  such that we have \begin{enumerate} 
            \item $\{ v, w \} = 0$ for any $v \in V^{\lambda}_{\mu, \C}:= V^{\lambda}_{\mu, \zee} \otimes_{\zee} \C$ and $w \in V^{\lambda}_{\mu', \C}$, for $\mu, \mu' \in \wts_{\lambda}$ with $\mu \neq \mu'$;
            \item $\{ v, w \} \in \zee$ for any $v, w \in V^{\lambda}_{\zee}$;
            \item $\{ v_{\lambda}, v_{\lambda} \} =1$ where $v_{\lambda}$ is the highest weight vector chosen in \eqref{Vzee:form};
            \item $\{ \xi_a v, w \} = \{ v, \xi_{-a} w \}$ for $v, w \in V^\lambda_{\C}$ and $a \in \rtsaffre$.  \end{enumerate}
For a uniqueness statement, we refer to \cite[Lemma 11.5]{kac:book}.

\subsection{Loop groups}
\label{subLoopGroups}
        
We now have all the ingredients necessary to introduce the precise notion of \textit{loop groups} used in this paper. Our construction follows the work of Garland \cite{gar:loopGroups}, which itself is a substantial generalization of Chevalley's construction in the finite-dimensional case (\textit{cf}. \cite{steinberg:book:lecturesChevalleyGroups} for an exposition of some of Chevalley's results) and uses in a crucial way the integral form $\ueg_{\zee}(\gext)$ constructed in \S \ref{subIntegralFormGarlandRep}. One of the nice features of Garland's construction is that from the Chevalley form, one can directly construct both an arithmetic subgroup as well as a `maximal compact' subgroup as we explain below.

\renewcommand{\gg}{\wh{G}}
In this section, let $k$ be any field of characteristic $0$. We will construct groups $\gg_k$ which, roughly speaking, are central extensions of the groups $G(k((t))$ if $G$ is the group corresponding to the Lie algebra $\mf{g}$. It will be important for us to actually construct a further extension of $\gg_k$ that we denote at $\gg_k^e$ by `exponentiating' the degree operator $\dd$. Said differently, the Lie algebra of $\gg_k$ corresponds to $\gext$ and that of $\gg_k^e$ to $\gext^e$. It is possible to make precise this analogy with the Lie algebras in some contexts but we do not need this here.

\tpoint{Exponentials and the polynomial loop group}

Fix $\lambda \in \Lambda_+$. Recall that we constructed a representation $V:=V^{\lambda}$ of the Lie algebra $\gext^e$ in \S \ref{subsubIrreducible} and introduced its integral form $V_{\zee}^{\lambda}$ in \S \ref{SubSubIntFormsRep}. In \S \ref{SubSubRepVk}, for any field $k$, we wrote $V_k:=V_{\zee} \otimes_{\zee} k$, which is a representation for the Lie algebra $\gext_k^e$. As noted above, the representation~$V$ is integrable (cf. \cite[Lem. 7.12]{gar:loopGroups}), \textit{i.e.}  for every vector $v \in V^{\lambda}$, there exists an integer $r > 0$ such that we have $\xi_{a_i}^r \cdot v  =  0$ for any simple root $a_i \in \Pi$, where we recall that $\xi_{\, a_i}$ is an element from our fixed Chevalley basis introduced in \S \ref{subsubChevalleyBasisLoop}, which lies in the root space $\mf{g}_{\zee}^{a_i}$. A similar statement also holds for $a_i$ replaced with any real root $a= \alpha + n \iota \in \rtsaffre$ with $\alpha \in \rts$ and $n \in \zee$, which we can see by conjugating $a$ to some $a_i$ using the Weyl group. For any such real root, this result allows us to define \be{chi:exp} \begin{array}{ccccc}\chi_{a} \left( s \right) & := & \exp \left( \xi_{a} s \right) & = & \sideset{}{}\sum\limits_{n \geqslant 0} \; \frac{\xi_{a}^n}{n!} s^n \end{array} \ee as an automorphism of $V_k,$ and, if we have $s \in \zee$, the corresponding element actually preserves $V_{\mathbb{Z}}$ as well. Using these operators we define an analogue (of the central extension) of $G(k[t, t^{-1}])$ by considering \be{}\label{incom:group} \begin{array}{lll} \gg_{k, \mathrm{pol}} & := & \la \, \chi_{a_i}(s) \; \mid \; i \in \lb 1, \ell+1 \rb, \; s \in k \, \ra, \end{array} \ee where $\left< \cdot \right>$ stands for `the group generated by.' We will deal with a completion of this group in which the ring of Laurent polynomials $k[t, t^{-1}]$ is replaced by the field of (formal) Laurent series $k((t))$.

\addtocounter{theorem}{1}

\tpoint{Exponentiation and Laurent series} 

To define the completed loop groups, we need the following extension of the fact that $V$ is integrable.

\begin{proposition}\cite[Lem. 7.16]{gar:loopGroups}
    Let $\alpha \in \rts$ be a classical root. For every vector $v \in V_{\zee}$, there exists an integer $n_0$ such that for every integer $n \geqslant n_0$, we have $\xi_{\alpha + n \iota}\,  \cdot v  =  0$.
\end{proposition}
\noindent As a consequence, for any Laurent series $\sigma(t) \in k (( t ))$, which we write as
    \begin{equation}
        \label{EqLaurentSeries}
            \begin{array}{ccc}
                    \sigma \left( t \right) & = & \sideset{}{}\sum\limits_{j \geqslant j_0} q_j t^j
                \end{array},
    \end{equation}
                
                \noindent where $j_0$ is an integer and $q_{j_0}$ is non-zero, possibly negative, and any classical root $\alpha \in \rts$, we may define
    \be{}\label{exp:comp} 
        \begin{array}{ccc}
            \chi_{\alpha} \left( \sigma \left( t \right) \right) & = & \sideset{}{}\prod\limits_{j \geqslant j_0} \; \chi_{\alpha + j \iota} \left( q_j \right),
        \end{array}
    \ee
\noindent as an element of $\mathrm{Aut}(V_k)$. The automorphisms in the right-hand side of \eqref{exp:comp} are pairwise commutative, following \cite[p. 47]{gar:loopGroups}.

\tpoint{The completed loop groups} \label{subsubCompletedLoopGroups} Having defined the automorphisms $\chi_{\alpha} \left( s \left( t \right) \right)$, we can construct the completed loop groups, following an approach similar to \eqref{incom:group}. We fix a dominant integral weight $\lambda \in \Lambda_+$.
            
\begin{definition}
    For any field $k$, we define the \textit{completed loop group} as \be{}\label{defCompletedLoopGroup} \begin{array}{lll} \widehat{G}_k & = & \left< \, \chi_{\alpha} \left( s \left( t \right) \right) \; \vert  \; \alpha \in \rts, \; \sigma \left( t \right) \in k (( t )) \, \right> \end{array}, \ee where once again the notation $\left< \cdot \right>$ means `the group generated by'.
\end{definition}
            
\noindent If an affine real root $a \in \rtsaffre$ is written as $\alpha + n \iota$, then note that we have \be{chi:a} \begin{array}{lll} \chi_a(s) & = & \chi_{\alpha}( s t^n). \end{array} \ee The group $\gg_k$ thus contains one-parameter subgroups $\{ \chi_a(s) \mid s \in k \}$ for each real root. We do not however have any corresponding one-parameter subgroup for the imaginary roots. Note also that the group $\gg_k$ above depends on the choice of representation, \textit{i.e.} on $\lambda$, but in a known (and rather mild) way (see \cite[Sec. 15]{gar:loopGroups}). Here, we fix some $\lambda$ such that $\la \lambda, \cc \ra > 0$, and drop it from our notation.

\tpoint{Remark on central extensions} In this formulation, the central extension is built in from the start, \textit{i.e.} we do not construct $\gg_k$ as central extension of some other group. However, one can see from \cite[\S 12]{gar:loopGroups} that these groups are actually central extensions of groups of the form $G(k((t)))$. In fact, if $G$ is simple, the groups introduced above are essentially the pushforward of a $K_2$ (or Steinberg--Matsumoto--Moore) universal extension of $G(k((t)))$ under the tame symbol \be{}\label{tameSymbol} \begin{array}{lll} k((t))^* \; \times \; k((t))^* & \longrightarrow & k^* \end{array} . \ee This fact does not play any role in our work here at the moment, but may be important in extending our results to an \textit{ind-pro} context, if one follows the analogy with the function field case.
            
\addtocounter{theorem}{1}

\tpoint{Extended groups} \label{subsubExtendedGroups} The full group analogue of $\gext^e_k$ is obtained from $\gg_k$ (or $\gg_{k,\mathrm{pol}}$) using a twist by a certain automorphism. For any $\tau \in k^*$, define $\eta \left( \tau \right) \in \aut(V_k)$ through its action on weight spaces:            
    \be{}\label{actionEtaTauWeightSpaces} 
        \begin{array}{llllll}
            \eta \left( \tau \right) \cdot v & = & \eta \left( \tau \right)^{\mu} v & = & \tau^{\mu \left( \mathbf{d} \right)} v & \text{ for every vector } v \in V_{k, \mu}.
        \end{array} \ee These automorphisms normalize $\gg_k \subset \aut(V_k)$ since one can check that we have
    \be{}\label{rotationLoopParameter}
        \begin{array}{llll}
            \eta \left( \tau \right) \chi_{\alpha} \left( \sigma \left( t \right) \right) \eta \left( \tau \right)^{-1} & = & \chi_{\alpha} \left( \sigma \left( \tau t \right) \right) & \text{ for } \sigma(t) \in k((t))
        \end{array} . \ee
    \begin{definition}
        \label{DefExtendedLoopGroup}
            The \textit{extended loop group} is defined by the semi-direct product
                \be{}\label{extendedLoopGroup}
                    \begin{array}{lll}
                        \widehat{G}^e_k & = & \widehat{G}_k \rtimes_{\eta} k^{\ast}
                    \end{array} , \ee the group $k^{\ast}$ acting on $\widehat{G}$ by $\eta \left( \tau \right)$. The elements of $\gg^e_k$ are written as $g \eta(\tau)$ with $g \in \gg_k$ and $\tau \in k^{\ast}.$ 
    \end{definition}
    A polynomial version $\gg_{k,\mathrm{pol}}^e$ of \eqref{extendedLoopGroup} is obtained by replacing $\sigma(t) \in k((t))$ with $\sigma(t) \in k[t, t^{-1}]$. Furthermore, we often fix $\tau$ and consider the following sub\emph{set}
        \be{}\label{gEtaTau}
            \begin{array}{lllllll}
                \gg^{\tau}_k & := & \widehat{G}_k \eta \left( \tau \right) & := & \left \{ g \eta \left( \tau \right), \; g \in \widehat{G}_k \right \} & \subset & \widehat{G}^e_k
            \end{array} , \ee with a similar definition for $\gg_{k,\mathrm{pol}} \eta \left( \tau \right)$.
            
            \newcommand{\hh}{\wh{H}} \newcommand{\uu}{\wh{U}}

\tpoint{Toral subgroup $\hh^e_k$}             
In this paragraph, let us present some important subgroups of $\gg^e_k$ which can be defined over any field. To do so, we first need to introduce certain elements in the group. For each affine root $a \in \rtsaffre$ and $s \in k^*$ we introduce the elements     
    \be{w:real}  
        \begin{array}{lllllll}
            w_{a} \left( s \right) & = & \chi_{a} \left( s \right) \chi_{- a} ( - s^{-1} ) \chi_{a} \left( s \right), & \text{ and } & h_{a} \left( s \right) & = & w_{a} \left( s \right) w_{a} \left( 1 \right)^{-1}.
        \end{array} \ee From \cite[Lemma 11.2.ii]{gar:loopGroups}, the elements $h_a(s)$ act as multiplication by $s^{ \la \mu, \av \ra}$ on the weight space $V_{k, \mu}$. In particular, we have $h_a ( s )^{-1} = h_a ( s^{-1} )$. Thus, with respect to an admissible basis of $V_k$, these are diagonal elements, and we define the abelian subgroup \be{torus} \begin{array}{lllll} \hh_k & := & \la \, h_a(s) \; \mid \; a \in \rtsaffre, \; s \in k^* \, \ra & \subset & \gg_k. \end{array} \ee Using \cite[Lem. 4.4]{gar:loopAlgebras}, one can argue that every element $h \in \hh_k$ can be written as \be{h:exp} \begin{array}{llll} h & = & \sideset{}{}\prod\limits_{i=1}^{\ell+1} \; h_{a_i}(s_i) ^{n_i}& \text{ for some scalars } \; s_i \in k^* \; \text{ and integers } \; n_1, \ldots, n_{\ell+1} \geqslant 0. \end{array} \ee This expression may not be unique and depends on which representation $V:= V^{\lambda}$ we pick to define the group. In favourable cases, \textit{e.g.} when $\lambda$ is picked so that the weight lattice spanned by $\wts_{\lambda}$ is equal to $\Lambda$, then the expression is unique (the \emph{simply-connected case}). We also define the extended torus \be{ext:torus} \begin{array}{lll} \hh_k^e & := & \hh_k \; \rtimes_{\eta} \; k^*. \end{array} \ee 
                    
\noindent For any elements $\mu \in \wts_{\lambda}$ and $x = h \eta(\tau) \in \hh_{k}^e$, with $h \in \hh_{k}$ and $\tau \in k^*$, we define \be{h:mu} \begin{array}{lllll} x^{\mu} & = & (h \eta(\tau))^{\mu} & := & \sideset{}{}\prod\limits_{i=1}^{\ell+1} \; s_i^{\la \mu, \av_i \ra}  \, \tau^{\la \mu, \dd \ra } \end{array} \ee with respect to \emph{any} decomposition (\ref{h:exp}). This expression is well-defined as one verifies from the fact that we have $x. v= x^{\mu} \, v$ for any vector $v \in V_{\mu}$. The definition of $x^{\mu}$ is extended to the integral span of $\wts_{\lambda}$.

\addtocounter{theorem}{1}

\tpoint{Torus in the real group} \label{baseFieldRealLog} Suppose we have $k = \mathbb{R}$. The subgroup $\widehat{H}_+$ of $\widehat{H}_{\R}$ is defined as \be{}\label{subgroupHPlusR} \begin{array}{lll} \widehat{H}_+ & = & \la \, h_{a_i} \left( s \right) \; \mid \; i \in \llbracket 1, \ell+1 \rrbracket, \; s > 0 \, \ra \end{array}. \ee The \textit{logarithm map} is induced by the equalities \be{}\label{logMap} \begin{array}{llll} \log h_{a_i} \left( s \right) & = & \left( \log s \right) \av_i & \text{ for } i \in \llbracket 1, \ell+1 \rrbracket \end{array} \ee which we extend as a group morphism to a map also denoted $\log: \widehat{H}_{+} \longrightarrow \mf{h}_{\mathbb{R}}$. We now set \be{}\label{extendedRealToral} \begin{array}{lll} \widehat{H}_{+}^e & = & \widehat{H}_{+} \; \rtimes \; \mathbb{R}^{*}_{+} \end{array}, \ee where $\mathbb{R}^{*}_{+}$ acts on $\widehat{H}_{+}$ by $\eta \left( \tau \right)$. Setting $\log \eta \left( \tau \right) = \left( \log \tau \right) \dd$ for any $\tau > 0$, we extend the logarithm map into a group morphism $\log: \widehat{H}_{+}^e \longrightarrow \mf{h}_{\mathbb{R}}^e$.

\addtocounter{theorem}{1}

\tpoint{$\bb_k$ and the (pro-)unipotent subgroup $\uu_k$}  The respective analogues of the upper and lower triangular subgroups (with respect to some fixed coherent basis of $V$) in $\gg_{k, \mathrm{pol}}$ are the groups generated by $\chi_{a}(s)$ for $a \in \rtsaff_{\mathrm{re}, \pm}$ and $s \in k$. We explain how to  complete the former (the upper triangular group). To do so, introduce the following elements of $\gg_k:$ for any $\alpha \in \rts$ and $\sigma \left( t \right) \in k ((t)) ^{\ast}$, \textit{i.e.} a non-zero element, we set
                  \be{def:w}  \begin{array}{lll}
                        w_{\alpha} \left( \sigma \left( t \right) \right) \; = \; \chi_{\alpha} \left( \sigma \left( t \right) \right) \chi_{- \alpha} ( - \sigma \left( t \right)^{-1} ) \chi_{\alpha} \left( \sigma \left( t \right) \right), & \text{ and } & h_{\alpha} \left( \sigma \left( t \right) \right) \; = \; w_{\alpha} \left( \sigma \left( t \right) \right) w_{\alpha} \left( 1 \right)^{-1}.
                    \end{array} \ee Note that the latter elements (when $\sigma(t)$ is not a scalar) \emph{do not} act diagonally. In fact, for any $\sigma(t) \in k[[t]]^*$ with $\sigma(0)=1$, we may use \cite[Lemma 8.1 and discussion\footnote{The result is written for $k=\R$ but works for a general field.} after]{gar:duke} to write \be{h:fact} \begin{array}{lll} h_{\alpha}(\sigma(t)) & = & \sideset{}{}\prod\limits_{n=1}^{\infty} h_{\alpha} ( 1 - b_n t^n) \end{array} \ee for uniquely determined scalars $b_n \in k^*$. The elements $h_{\alpha} ( 1 - b_n t^n)$ now act on $V$ as (parts of) certain vertex operators\footnote{We owe this observation to H. Garland and it is recorded in some form in \cite[\S 4.7]{bgkp}.}.
                    
    \begin{definition}
        The subgroup $\uu_k \subset \gg_k$ is defined as the subgroup generated by the following elements: 
            \begin{itemize}
                \item[$\bullet$] $\chi_{\alpha} \left( \sigma \left( t \right) \right)$ with $\alpha \in \rts+$ and $\sigma(t) \in k[[t]]$;
                    
                \item[$\bullet$] $\chi_{\alpha} \left( \sigma \left( t \right) \right)$ with $\alpha \in \rts_-$ and $\sigma(t) \in \; t \, k[[t]]$;
            \end{itemize}
    \end{definition}
\noindent The group $\hh_k$ normalizes $\uu_k$, so we may define $\bb_k:= \uu_k \rtimes \hh_k$ and $\bb^e_k:= \uu_k \rtimes \hh^e_k$, which play the role of a Borel subgroup for $\gg_k$. In fact, as explained in \cite[Sec. 12.1]{gar:loopGroups}, with respect to an integral admissible coherently ordered basis (see \S \ref{subsubAdmissible}), the elements of $\uu_k$ can be seen as (infinite) upper-triangular matrices with only ones on the diagonal and the elements of $\bb_k$ or $\bb_k^e$ can be seen as upper-triangular matrices.

\tpoint{Iwahori--Matsumoto type factorization} We will also need the factorization \cite[Sec. 18.11]{gar:loopGroups} for elements from $\uu_k$ reminiscent of the Iwahori--Matsumoto factorization for the pro-unipotent radical of an Iwahori subgroup in the realm of $p$-adic groups. This allows us to put ``coordinates'' on $\uu_k$.
    \begin{proposition}
        Every element $x \in \uu_k$ has a factorization 
            \be{}\label{Iwahori-MatsumotoFactorization} 
                \begin{array}{lll}
                        x & = & \left( \sideset{}{}\prod\limits_{\alpha \in \rts_+} \chi_{\alpha} \left( \sigma_{\alpha} \left( t \right) \right) \right) \, \left( \sideset{}{}\prod\limits_{j=1}^{\ell} \, h_{\alpha_i}(\sigma_i(t)) \, \right)  \left( \sideset{}{}\prod\limits_{\alpha \in \rts_+} \chi_{-\alpha} \left( \sigma'_{\alpha} \left( t \right) \right) \right)
                \end{array} ,
            \ee where, for any $\alpha \in \rts_+$, we have $\sigma_{\alpha}(t) \in k[[t]]$ and $\sigma'_{\alpha}(t) \in t k[[t]]$, as well as $\sigma_i(t) \in k[[t]]^*$ with $\sigma_i(0) =1$ for any $i \in \llbracket 1, \ell \rrbracket$. Moreover, the elements $\sigma_{\alpha}, \sigma'_{\alpha}, \sigma_i$ are uniquely determined once we fix an order on $\rts_+.$   
    \end{proposition} 
    
\tpoint{$\bb^-_k$ and its unipotent subgroup $\uu^-_k$}  We can also define the analogue of lower triangular subgroups. Note that the definition is simpler here as no completions are involved; the structure of these groups is however much more mysterious. 

    \begin{definition} \label{def:NegUnip}
        The subgroup $\uu^-_k \subset \gg_k$ is defined as being generated by 
          $\chi_{a}(s)$ for $a \in \rtsaff_{re, -}$ and $s \in k.$ We also set $\bb_k^{-}:= \uu^-_k \rtimes \hh_k $ and $\bb_k^{-, \, e}:= \uu^-_k \rtimes \hh_k^e$.
    \end{definition}

\noindent In \S \ref{subsubNegCompletions} we explain how to complete the groups $\uu^-_k$ in some sense. Though these no longer lie in $\gg_k$ or even act on~$V_k$, one can still use the completions for certain purposes.

\newcommand{\kk}{\wh{K}}

\tpoint{Maximal compact subgroup} \label{subsubMaximalCompact} Assuming we have $k = \mathbb{R}$, there is another subgroup of $\gg_{\mathbb{R}}$ we must define, which plays the role of a maximal compact subgroup. We set \be{K:inf} \begin{array}{lll} \kk & :=  & \{ \, k \in \gg_{\R} \; \mid \; \{  v , k w \} = \{ k^{-1} v, w \} \; \text{ for all } \; v, w \in V_{\mathbb{R}} \}. \end{array} \ee One can check, as in \cite[Lemma 11.2(i)]{gar:loopGroups}, that we have $w_{a}(\pm 1)\in \kk$ and also \be{kom:unip} \begin{array}{lll} \kk \; \cap \; \uu_{\R} & = & 1. \end{array} \ee

\addtocounter{theorem}{1} 

\newcommand{\Gam}{\wh{\Gamma}}

\tpoint{Arithmetic subgroup $\Gam$ } \label{subsubArithmeticSubgroup}Finally, let us introduce the last subgroup of $\gg_{\R}$ we need. We set \be{def:Gamma} \begin{array}{lll} \Gam & := & \{ \, \gamma \in \gg_{\R} \; \mid \; \gamma \, V_{\zee} = V_{\zee} \, \}. \end{array} \ee Note that this group also plays a role in defining the \textit{Siegel sets}. For more information, the reader is referred to \cite[Sec. 19]{gar:loopGroups}. We will not need this notion here.

         \addtocounter{theorem}{1}           
   \newcommand{\pol}{\mathrm{pol}}
   
\tpoint{Iwasawa decomposition}\label{IwasawaDecomp} In what follows, we will need a way to conveniently write the elements of $\gg_{\mathbb{R}}$, known as the \textit{Iwasawa decomposition}. It is an infinite-dimensional generalization of the Graham--Schmidt procedure from linear algebra.

\begin{theorem}[Iwasawa decomposition] \cite[Lem. 16.14]{gar:loopGroups}
    We have $\gg_{\R} = \kk \bb_{\R}$, and moreover every element $g \in \gg_{\R}$ can be uniquely factored as \be{}\label{iwasawaFactorization} \begin{array}{llll} g & = & k \, u \, h & \text{ for some elements } k \in \kk, \; u \in \uu_{\R}, \; h \in \hh_+ \end{array}. \ee
\end{theorem}

\noindent It should be noted, with respect to \cite[Lem. 16.14]{gar:loopGroups}, that the order of $\uu_{\R}$ and $\hh_+$ can be reversed, since $\hh_+$ normalizes $\uu_{\R}.$ We can also state a version of this decomposition for the extended groups of \S \ref{subsubExtendedGroups}. For a fixed $\tau \in \R^*$, every $x = g \eta(\tau) \in \gg^{\tau}_{\R}$ can be written uniquely as \be{extended-iwasawa} \begin{array}{lll} g \eta(\tau) & = & k \; u \; h \; \eta(\tau) \end{array} \ee with $k \in \kk$, $h \in \hh_+,$ and $u \in \uu_\R$. We could also exchange the order of the component in $\hh_+ \eta(\tau)$ with the component in $\uu_k$ since the extended torus $\hh^e_{\R}$ also normalizes $\uu_{\R}$. Let us also remark here that, using the Iwasawa decomposition, one can show (in analogy with \eqref{B:G}) that we have \be{bunG:B} \begin{array}{lll} \kk \setminus \gg^{\tau}_{\R}  / \Gamma \cap \bb_{\R} & = & (\kk \cap \bb_{\R}) \setminus \bb^{\tau}_{\R}  / \Gamma \cap \bb_{\R}, \end{array} \ee where we have adopted the convention \be{B:tau} \begin{array}{lll} \bb^{\tau}_{\R} & := & \bb_{\R}  \eta(\tau). \end{array} \ee

\tpoint{Remark on the Iwasawa decompositions with respect to $\bb^-_{\R}$} \label{subsubNegIwasawa} It is not true however that every element in $x \in \gg_{\R}$ (or $\gg^{\tau}_{\R})$ can be factored as $x = k b^-$ with $k \in \kk$ and $b^- \in \bb^-_{\R}$. Indeed, this stems from the following issue: the Iwasawa decomposition is essentially a consequence of a rank 1 computation and a certain Bruhat decomposition (see \cite[\S16]{gar:loopGroups} or \cite[Ch. 8]{steinberg:book:lecturesChevalleyGroups}). In order to have an Iwasawa decomposition with respect to $\bb^-_{\R}$ one needs a Bruhat decomposition of $\gg_{\R}$ into $(\bb^-_{\R}, \bb^-_{\R})$-double cosets, which fails in the completed group. However, if one replaces $\gg_{\R}$ with the more `symmetric' object $\gg_{\pol, \R}$ (which is more problematic from the point of view of reduction theory), then Bruhat decompositions hold with respect to both $\bb^-_{\R}$ and $\bb_{\R}$, and we have Iwasawa decompositions with respect to both $\bb_{\R}$ and $\bb_\R^-$. If we write, as above, $\bb^{-, \tau}_{\R}:= \bb^{-}_{\R}\eta(\tau)$, our main finiteness result (see Theorem \ref{thm:main-fin} below) holds for elements in \be{neg:bunB}  (\kk \cap \bb^-_{\R}) \setminus \bb^{-, \tau}_{\R}  / \Gamma \cap \bb^-_{\R}, \ee which is not equal to $\kk \setminus \gg^{\tau}_{\R}  / \Gamma \cap \bb^-_{\R}$. Rather, one has \be{}\label{} \begin{array}{lll} (\kk \cap  \gg^{\tau}_{\pol, \R} ) \setminus  \gg^{\tau}_{\pol, \R}  / \Gamma \cap \bb^-_{\R} & = & (\kk \cap \bb^-_{\R}) \setminus \bb^{-, \tau}_{\R}  / \Gamma \cap \bb^-_{\R}. \end{array} \ee

  \newcommand{\trr}{\twoheadleftarrow}
  \newcommand{\bun}{\mathrm{Bun}}
  \newcommand{\hbun}{\overline{\mathrm{Euc}}}
  
\newcommand{\naff}{\wh{\mf{n}}}

\section{Pro-Hermitian bundles from loop groups} \label{SecProHermitian}

In this section, we will use the results from sections \ref{SecHermitianLattices} and \ref{SecLoopStuff} to show how to construct a pro-Hermitian vector bundle over $\Spec \, \mathbb{Z}$ from an element in the loop group $\gg_{\mathbb{R}}$. Our construction factors to give a map \be{pro-quot} \begin{array}{lllll} \mc{X}_{\bb} & := & \kk \setminus \gg_{\R} \eta(\tau) / \Gam \cap \bb^-_{\mathbb{R}} & \stackrel{\Psi}{\longrightarrow} & \mathrm{iso}( \provectZ) \end{array} \ee where the right-hand side is the set of isomorphism classes of pro-Hermitian vector bundles over $\Spec \, \zee.$ The map $\Psi$ depends on a choice of a representation $V$ of the group $\bb^-$, and we are especially interested in the case when $V$ is a highest weight representation (\textit{cf.} \S \ref{subGarlandRep}), which are representations of $\gg^e$ viewed as representations of $\bb$. Assuming we have $\left \vert \tau \right \vert<1,$ we can additionally show that the pro-Hermitian bundles we construct are theta-finite (\textit{cf.} Theorem \ref{thm:main-fin} ). We remark that the natural setting for our construction is actually the adelic one, but since we are only considering the number field $K=\Q$ here, a version of strong approximation for loop groups (\textit{cf}. \cite[Appendix B]{garland:article:loopGroups2HilbertModular}) allows to just work over the infinite place, \textit{i.e.} with real groups. The adelic picture will be presented in a future work.

\subsection{Constructing pro-Hermitian bundles}
Let $\lambda \in \Lambda_+$ be a dominant integral weight with $\lambda \left( \cc \right) > 0$, to which we associate the highest weight representation $V$ of $\gext^e$ constructed in \S \ref{subGarlandRep}. We pick an integral admissible coherently ordered basis $\mc{B}$ (cf. \S \ref{subsubAdmissible}) of $V$, write $V_{\zee}$ for its $\zee$-span, and set $V_{\mathbb{R}} := V_{\zee} \otimes_{\zee} \R$. In this section, we work with the completed loop group $\gg_{\R} \subset \aut(V_{\R})$ defined in \S \ref{subsubCompletedLoopGroups}. We defined in \S \ref{subsubArithmeticSubgroup} the arithmetic subgroup $\Gam \subset \gg_{\R}$ preserving the lattice $V_{\zee}$. In \S \ref{subsubHermitianStructure}, we defined a Hermitian inner product on $V = V_{\mathbb{R}} \otimes_{\mathbb{Z}} \mathbb{C}$, and, in \S \ref{subsubMaximalCompact}, denoted by $\wh{K} \subset \gg_{\mathbb{R}}$ the subgroup of isometric automorphisms of $V$. We write $\left \Vert \cdot \right \Vert$ for the norm associated to the inner product. Finally, we consider a scalar $0  < \tau < 1$, and consider the subset $\gg_{\R} \eta(\tau) \subset \gg^e_{\R}$. The dominant weight $\lambda$ being fixed, we drop it from the notations.

\tpoint{Level quotients} \label{subsubLevelQuotients}  Recall from \S \ref{subsubDepthLevel} that we denote by $\wts(V)$ the set of weights of $V$. We have also defined the level of weight $\mu \in \wts(V)$ in \eqref{depth:level}. We now define the subspaces $V_n \subset V$ and $V^n \subset V$ as \be{level-subspaces} \begin{array}{lllllll} V_n & := & \bigoplus_{\substack{\mu \in \wts(V), \\ \, \lev(\mu) \leqslant n}} \, V_{\mu} & \; \text{ and } \; & V^n & := & \bigoplus_{\substack{\mu \in \wts(V), \\  \lev(\mu) > n}} \, V_{\mu}. \end{array} \ee From the basic properties of $\{ \cdot, \cdot \}$ described in \S \ref{subsubHermitianStructure}, the space $V_n$ and $V^n$ are orthogonal for $\{ \cdot, \cdot \}$. Note that we have the following identification \be{quot:sub-id} \begin{array}{lll} V_n & \cong & V / V^n. \end{array} \ee  

\noindent From the definition of $V_{\mathbb{Z}}$ as the $\mathbb{Z}$-span of our integral admissible coherently ordered basis $\mc{B}$, we find that \be{Vn:zee} \begin{array}{lllllll} V_{n, \zee} & := & V_n \cap V_{\zee} & \; \text{ and } \; & V^n_{\zee} & := & V^n \cap V_{\zee} \end{array} \ee are the $\zee$-span of elements from $\mc{B}$ which have level lower than $n$, or strictly greater than $n$, respectively. Moreover,  the identification (\ref{quot:sub-id}) restricted to integral lattices yields an identification \be{quot:sub-id:int} \begin{array}{lll} V_{n, \zee} & \cong & V_{\zee}/V^{n}_{\zee}. \end{array} \ee Finally let us define the subspaces and sublattices \be{V[n]} \begin{array}{lllllll} V[n] & := &  \bigoplus_{\substack{\mu \in \wts(V), \\ \, \lev(\mu) = n}} & \; \text{ and } \; & V_{\zee}[n] & = & V[n] \cap V_{\zee}. \end{array} \ee Note that the kernel of the natural surjection \be{qn} \begin{array}{llllll} q_{n} & : & V / V^{n+1} & \longrightarrow & V/V^{n} & \text{ for } \; n \geq 0 \end{array} \ee is isomorphic, under \eqref{quot:sub-id}, to the subspace $V[n+1].$ Moreover, we also have \be{VZ:ker} \begin{array}{llll} V_{\zee}[n+1] & = & \ker\left( q_{n, \zee} : V_{\zee} / V^{n+1}_{\zee} \rr V_{\zee} / V^{n}_{\zee} \right) & \text{ for } n \geq 0 . \end{array} \ee We shall often write \be{Y:def} \begin{array}{lllllll} Y_n & := & V / V^n & \; \text{ and } \; & Y_{n, \zee} & := & V_{\zee} / V^n_{\zee} \end{array}. \ee

\addtocounter{theorem}{1}

\tpoint{A completion of $\uu^-_{\R}$} \label{subsubNegCompletions} In Definition \ref{def:NegUnip}, we defined a group $\uu^-_{\R}$ acting on $V.$ In \cite[\S 4]{bgkp}, following the construction of the Iwahori--Matsumoto coordinates from \cite[\S 18]{gar:loopGroups} on $\uu^-_{\R},$ a completion $\uu_{\R, \comp}^-$ was introduced together with a map \be{}\label{} \begin{array}{lll} \uu^-_{\R} & \longrightarrow & \uu_{\R, \comp}^-. \end{array} \ee  One knows that we have $\uu^-_{\R} (V^n) \subset V^n$, so $\uu^-_{\R}$ induces an action on $V/ V^n$. Although the completion $\uu_{\R, \comp}^-$ does not act on $V$, as applying a negative unipotent from this completion can yield non-zero components in infinitely many different weight spaces, it still acts on the finite-dimensional quotients $V/ V^n$ for each $n \geq 0,$ and in fact for each such $n$ we have a diagram \be{um:act} \begin{tikzcd} \uu^-_{\R} \ar[rr] \ar[dr] & & \uu_{\R, \comp}^- \ar[dl] \\ & \aut(V/V^n) &  \end{tikzcd}. \ee As stated in \cite[Thm 4.4]{bgkp}, every element of the completion, and hence every element $x \in \uu^-_{\R}$, has a factorization in $\uu^-_{\R, \comp}$ of the following form \be{uminus:factorization} \begin{array}{lll} x & = & \sideset{}{}\prod\limits_{k \geq 0} u^-(k), \end{array} \ee where each piece $u^-(k)$ is uniquely determined from $x$ and acts on $V/ V^n$ as \be{}\label{uminus:act} \begin{array}{llll} u^-(k) w & = & w + w' & \text{ if } \; w \in V[j] \; \mathrm{ mod } \; V^n \; \text{ and where } \;  w' \in V[j+k] \; \mathrm{ mod } \; V^n. \end{array} \ee The proof is analogous to the corresponding result for $\uu_k$ and we refer to \cite[\S 18, pp. 96-97]{gar:loopGroups} for more details. Note that by construction $u^-(0)$ consists of a product \be{}\label{uminus:0} \begin{array}{ll} \sideset{}{}\prod\limits_{\alpha \in \rts_+ } \chi_{\alpha}(s_{\alpha}) & \text{ for } \; s_{\alpha} \in \R. \end{array} \ee If we fix an ordering on $\rts_+$, this description can even be made unique.    
\addtocounter{theorem}{1}

\tpoint{Quotient norms and a projective system} \label{subsubConstructProSystem} For any $x \in \gg_{\R}^{\tau}$ we define the \textit{twisted norm} on $V$ as \be{twist:norm} \begin{array}{llll} \left \Vert v \right \Vert_x & := & \left \Vert x v \right \Vert & \text{ for } v \in V. \end{array} \ee We can equip each quotient vector space $Y_n:=V / V^n$ with the quotient norm (\textit{cf}. \cite[Sect. 1.1.3]{chen-moriwaki:book:arakelov}), which we write as $\left \Vert \cdot \right \Vert_{V/V^n, \, x}.$ By definition, we have \be{def:quot} \begin{array}{lllll} \left \Vert v \right \Vert_{V/V^n ,\,  x} & := & \inf\limits_{y^n \in Y^n } \left \Vert v + y^n \right \Vert_x & = & \inf\limits_{y^n \in Y^n } \left \Vert x( v + y^n) \right \Vert, \end{array} \ee and we further set \be{Yn:euc} \begin{array}{lllll} \ov{Y}_{n,x} & := & (Y_{n, \zee}, \left \Vert \cdot \right \Vert_{V/V^n, x} )  & \in & \vectZ \end{array} \ee  When the norm is not twisted, \textit{i.e.} when we have $x=1$, we write $\left \Vert \cdot \right \Vert_{V/V^n}:= \left \Vert \cdot \right \Vert_{V/V^n, 1}.$ Using the general properties of quotient norms (\textit{cf.} \cite[Proposition 1.1.4]{chen-moriwaki:book:arakelov}), it follows that under the natural surjection \be{}\label{surjQuotient} \begin{array}{llll} V / V^n & \longrightarrow & V / V^{m} & \text{ with } n \geqslant m, \end{array} \ee the norm $\left \Vert \cdot \right \Vert_{V/V^m, x}$ is the quotient norm obtained from $\left \Vert \cdot \right \Vert_{V/V^n, x}$ under the quotient map above. Thus the elements $\ov{Y}_{n, x}$ fit into a \emph{projective} system \be{}\label{proj:sys} \begin{array}{lllllllll} \ov{Y}_{\bullet,  x} & : & \ov{Y}_{0, x} & \stackrel{q_0}{\longleftarrow} & \ov{Y}_{1, x} & \stackrel{q_1}{\longleftarrow} & \ov{Y}_{2, x} & \longleftarrow & \cdots, \end{array} \ee of admissible surjective morphisms (\textit{cf}.\S  \eqref{admissibleProjSyst}) and hence defines an element of $\provectZ$ (\textit{cf}. \S \ref{proHermitianFromProjSyst}). From an element $x \in \gg_{\R}^{\tau}$ and our fixed highest weight representation $V$, we thus constructed an element \be{}\label{def:Psi} \begin{array}{lllll} \Psi \left( x \right) & = & \overline{Y}_{\bullet,  x} & \in & \provectZ. \end{array} \ee

\addtocounter{theorem}{1}

\tpoint{ The $(\kk, \Gam \cap B^-)$-invariance of $\Psi$} \label{subsubGammaKInvariance} We will now investigate the natural domain and range of the map $\Psi$ constructed above in \eqref{def:Psi}. First note that for any element $k \in \kk$, having \be{k-inv} \begin{array}{llll} \left \Vert v \right \Vert_{kx} & = & \left \Vert v \right \Vert_x & \text{ for any } v \in V, \end{array} \ee implies that we have \be{Psi:K-inv} \begin{array}{lll} \Psi (kx ) & = & \Psi(x). \end{array} \ee 

\noindent Now, consider an element $\gamma \in \Gam \cap \bb^-.$ For any element $b \in \bb^-$, we have $b V^n = V^n$, so we also have \be{}\label{} \begin{array}{llll} \min\limits_{y^n \in Y^n} \left \Vert x (\gamma v + y^n) \right \Vert  & = & \min\limits_{y^n \in Y^n} \left \Vert x \gamma (v + y^n) \right \Vert & \text{ for any } v \in V. \end{array} \ee In other words we have \be{quotnorm:gam} \begin{array}{llll} \left \Vert \gamma v \right \Vert_{V/V^n, x} & = & \left \Vert v \right \Vert_{V/V^n, x\gamma} & \text{ for any } v \in V . \end{array} \ee

\noindent The element $\gamma \in \Gam \cap \bb^-$ thus induces an (isometric) isomorphism of metrized lattices \be{} \begin{array}{lllll} m_{\gamma} & : & V_{\zee} / V^n_{ \zee} & \longrightarrow & V_{\zee} / V^n_{\zee}. \end{array} \ee Furthermore, these maps $m_{\gamma}$ commute with the maps $q_n$, and therefore induce an isometry of projective systems $\overline{Y}_{\bullet, \gamma x} \simeq \overline{Y}_{\bullet, x}$, which translates into an isometry of the associated pro-Hermitian vector bundles.

\begin{proposition} 
\label{prop:PsiKGamma} 
    For any fixed highest weight representation $V$ of $\gg_{\R}^{\tau}$, the map $\Psi$ from \eqref{def:Psi} descends to a map, denoted by the same name $\Psi:= \Psi_V$ \be{Psi:factored} \begin{array}{lllll} \mc{X}_{\bb^-} & := & \kk \setminus \gg^{\tau}_{\R} / \Gam \cap \bb^- & \stackrel{\Psi_V}{\longrightarrow} & \mathrm{iso} \, \left( \provectZ \right), \end{array} \ee where the right-hand side denotes the set of isomorphism classes in $\provectZ$. 
\end{proposition}

\subsection{Theta-finiteness} 

The main result of this paper, is that the pro-Hermitian vector bundle $\Psi \left( x \right)$ is theta-finite for any element $x \in \bb^{\tau, -}$, where $\tau$ satisfies $0 < \tau < 1 $. To that effect, we will verify that the strong summability condition \eqref{strongSummabilityCondition} is satisfied for the projective system \eqref{proj:sys} for every $\varepsilon > 0$. Let us first make several simplifications.

\tpoint{Removing the quotient norms} We begin by writing an Iwasawa decomposition (\textit{cf.} \S \ref{IwasawaDecomp}) \be{}\label{iwa:neg} \begin{array}{llll} x & = & u^- \, h \, \eta(\tau) & \text{ where }  h \in \hh_+, \text{ and } u^- \in \uu^-. \end{array} \ee To check the strong summability condition, we need to consider the kernels of the maps $q_n$ from the system \eqref{proj:sys}, and their invariant $h_{\vartheta}^0$. Using \eqref{VZ:ker}, we have \be{}\label{} \begin{array}{lll} \overline{\ker (q_n)} & = & (V \left[ n \right], \; \left \Vert \cdot \right \Vert_{V/V^n, x}) \end{array} \ee where $\left \Vert \cdot \right \Vert_{V/V^n, x}$ is the quotient metric on $V \left[ n \right]$ when identified to $V / V^{n}$. We continue to denote by $\left \Vert \cdot \right \Vert_{x}$ the restriction of the metric on $V$ twisted by $x$ to the subspace $V \left[ n \right]$.

\begin{lemma} \label{lem:quot-norm} 

    \begin{enumerate}  
        \item If $x \in \wh{H}^+ \eta(\tau)$ and $v \in V_n$, we have \be{}\label{h:quot} \begin{array}{lll} \left \Vert v \right \Vert_{V/V^n, x} & = & \left \Vert v \right \Vert_x. \end{array} \ee

        \item Assume $u \in U^-$ has a factorization as $u^- = \prod_{n \geq 0} u^-(n)$ as in (\ref{uminus:factorization}). Then we have \be{}\label{unip:no-quot} \begin{array}{llll} \left \Vert u^- v \right \Vert_{V/V^n} & = & \left \Vert u^-(0) v \right \Vert & \text{ for any } v \in V \left[ n \right] \; = \; \ker \, q_n \end{array} \ee where we recall that $u^-(0)$ is an automorphism of $V \left[ n \right]$. 
    \end{enumerate} 

\end{lemma}

\begin{proof}

    First, we note that for any element $v \in V_n$, we have $\left \Vert v \right \Vert_{V/V^n,x} = \left \Vert v \right \Vert_x$, since $V_n$ is orthogonal to $V^n$ for the metric twisted by the diagonal element $x$, thereby giving the first part of the lemma. Let us move on to the second. Using \eqref{uminus:act} we have $u^- v = u^-(0) v \mod V^n$, which gives \be{}\label{inequalityUnipotent} \begin{array}{lllllll} \left \Vert u^- v \right \Vert_{V/V^n} & = & \inf\limits_{w \in V^n} \left \Vert u^- v + w \right \Vert & = & \inf\limits_{w \in V^n} \left \Vert u^-(0) v + w \right \Vert & = & \left \Vert u^-(0) v \right \Vert \end{array} \ee since $u^-(0)$ preserves $V \left[ n \right]$, which is orthogonal to $V^n$. The proof of the lemma is thus complete.

\end{proof}

\noindent In the study of the projective system \eqref{proj:sys}, with the norms twisted by $x$, we can use \eqref{unip:no-quot} to make a simplification, by replacing the full negative unipotent component of $x$ from \eqref{iwa:neg} with $u^-(0)$, and by removing the quotient norms.

\tpoint{Diagonal elements} Having taken care of the quotient norms and most of the unipotent part of~$x$ from \eqref{iwa:neg} in Lemma \ref{lem:quot-norm}, we now assume that we have $x = h \eta(\tau) u^-(0)$, with $u^-(0)$ as in \eqref{uminus:factorization}.

\begin{lemma} 
    \label{lem:expandH} 
    Consider an element $h \in \wh{H}^+$, a scalar $\tau$ with $0 < \tau < 1$, and write \be{}\label{} \begin{array}{lllll} \log(h \eta(\tau) ) & = & H + \left( \log \tau \right) \dd & \in & \haff^e_{\R}, \end{array} \ee with $H \in \haff$, as in \S \ref{baseFieldRealLog}. There exist two positive constants $C_1 := C_1 \left( \rts, x, \lambda \right) > 0$ and $A_1 := A_1 \left( \rts, x, \lambda \right) > 0$ such that, for any vector $v \in V \left[ n \right]$, we have \be{}\label{hv:lower-bound} \begin{array}{llll} \left \Vert h \eta \left( \tau \right) \, v \right \Vert & \geqslant & C_1 \tau^{- n} e^{- A_1 \sqrt{n}} \left \Vert v \right \Vert \end{array} . \ee
\end{lemma}

\begin{proof}
    We note that any weight $\mu \in \wts(V)$ of level $n$ may be written as $\mu = p \Lambda_{\ell+1} + \overline{\mu} - n \delta$, as in (\ref{mu:form}). For any such weight, and any vector $v \in V_{\mu}$, we have \be{}\label{} \begin{array}{lllllll} xv & = & (h \eta \left( \tau \right))^{\mu} \, v  & = &  \tau^{-n} h^{\mu} \, v & = & \tau^{-n} e^{p m} \, e^{\la \overline{\mu}, H_o \ra  } \, v, \end{array} \ee where we have written $H \in \haff= \mf{h} \oplus \R \cc$ as \be{H:dec} \begin{array}{llll} H & = & H_o + m \, \cc & \text{ with } \; m \in \R, \; H_o \in \mf{h}. \end{array} \ee The linear-quadratic inequality \eqref{lin-quad} and the definition of the operator norm on $\mf{h}^{\ast}$ now yield \be{cs:lq} \begin{array}{lllll} \left \vert \la \overline{\mu}, H_o \ra \right \vert & \leqslant & \left \vert \overline{\mu} \right \vert \, \left \vert H_o \right \vert & \leqslant & A_1 \sqrt{n} \end{array} \ee with a constant $A_1:= A_1(\rts, x, \lambda)>0.$ The result follows from this, the last constant being $C_1 := e^{pm}$.
\end{proof}

In the study of the projective system \eqref{proj:sys}, the norms being twisted by $x = h \eta(\tau) u^-(0)$, with $u^-(0)$ as in \eqref{uminus:factorization}, Lemma \ref{lem:expandH} will allow us to remove the diagonal component $h$ of $x$ by finding upper-bounds of the finite-rank invariants $h_{\vartheta}^0$.

\tpoint{Remaining unipotent component} Next, we want to establish a result which asserts that the $u^-(0)$ part of unipotent elements cannot shrink the norm of an element too much. 

\begin{lemma} 
    Consider an element $u^-(0) = \prod_{\alpha \in \rts_-} \chi_{\alpha}(s_{\alpha})$ as in \eqref{uminus:0}, the product being with respect to a fixed order. There exist constants $C_2:= C_2(x, \rtsaff, \lambda)>0$ and $A_2:= A_2(x, \rtsaff, \lambda) > 0$ such that we have  \be{}\label{unip:lowerbound} \begin{array}{lll} \left \Vert u^-(0) v \right \Vert & \geqslant & C_2 e^{- A_2 \sqrt{n}} \left \Vert v \right \Vert \end{array}. \ee
\end{lemma} 

\begin{proof} 
    Let us write the argument for the case where $u^-(0)$ is replaced by a positive unipotent of the form \be{}\label{posUnipotent} \begin{array}{lll} u & = & \sideset{}{}\prod\limits_{\alpha \in \rts^+} \chi_\alpha(s_{\alpha}) \end{array} \ee to align more closely with the available literature. The same argument works for $u^-(0)$. To begin, let us assume our element is given by a single $\chi_{\alpha}(s)$ for some real number $s > 0$. Then from \cite[(3.1) p.214]{gar:duke}, there exists a scalar $s' > 0$ (explicitly computed from $s$, see \emph{loc. cit.}) so that we may write a polar decomposition \be{}\label{polarDecomp} \begin{array}{llll} \chi_{-\alpha}(s) \chi_\alpha(s) & = & k^{-1}  h_{\alpha}(s') k & \text{ where } k \in \kk \text{ stabilizes } V[n]. \end{array} \ee Hence, we have \be{}\label{rank1:norm} \begin{array}{lllll} \left \Vert \chi_{\alpha}(s) v \right \Vert^2 & = & \la \chi_{-\alpha}(s) \chi_\alpha(s) \, v, \; v \ra & = & \la h_{\alpha}(s') k \, v, \; k \, v \ra \end{array}. \ee As $k$ preserves norms on $V[n]$, we need only estimate how $h_{\alpha}(s')$ affects norms on $V[n]$. To that effect, we can use the same argument as in the proof of the previous Lemma (\textit{i.e.} the argument based on the linear-quadratic inequality). The case of a (finite) product of elements follows from by induction.

\end{proof}

\tpoint{Rank estimate} Finally, we need to find an estimate on the dimension of $V[n].$ 

\begin{lemma}
    There exist constants $A_3:= A_3(\rtsaff, \lambda) > 0$ and $C_3:=C_3(\rtsaff, \lambda) > 0 $ such that we have \be{dim:V[n]} \begin{array}{lll} \dim_{\C} V[n] & \leqslant & C_3 e^{A_3 \sqrt{n}}. \end{array} \ee 
\end{lemma} 

\begin{proof} Writing $\wts_{\lambda}[n]$ for the set of weights of level $n$ in the representation $V$, we have \be{d:mult} \begin{array}[t]{lllll} \dim \, V[n] & = & \sideset{}{}\sum\limits_{\mu \in \wts[n]} \mult(\mu) & \leqslant & \left \vert \wts_{\lambda}[n] \right \vert \; \; \max\limits_{\mu \in \wts_{\lambda}[n]} \; \mult(\mu), \end{array} \ee where $\mult(\mu)$ is the multiplicity of the weight $\mu$, defined in \eqref{mult:wt} as the dimension of the corresponding weight space. From the Kac--Peterson estimate Theorem \ref{thm:KPbound} for these weight-multiplicities, it suffices to show that the number of distinct weights of level $n$ grows as a polynomial in $n$.  Indeed, writing \be{}\label{} \begin{array}{lll} \mu & = & p \Lambda_{\ell+1} + \overline{\mu} - n \delta \end{array} \ee as in \eqref{mu:form}, we know that $\overline{\mu}$ is an integral linear combination of $a_1, \ldots, a_{\ell}$, and, having $\mu \preccurlyeq \lambda$, we must have \be{}\label{fundamentalWeightCentralElement} \begin{array}{lllll} p & = & \la \mu, \cc \ra  & = & \la \lambda, \cc \ra, \end{array} \ee which then depends only on the highest weight $\lambda$. Now (\ref{pol-est}) states that we have \be{}\label{pol-est:2} \begin{array}{lll} \left \vert \overline{\mu} \right \vert^2  & \leqslant &  2  \, p \, n + \left \vert \lambda \right \vert^2. \end{array} \ee The vector space $\mf{h}^{\ast}$ being finite dimensional of dimension $\ell$ independent from $n$, the norms on $\mf{h}^{\ast}$ are equivalent, with implicit constants independent from $n$. In particular, the normalized Killing form $(\cdot, \cdot)$ restricted to $\mf{h}^{\ast}$ is positive-definite, thus induces a norm which is equivalent to the sup norm of the coefficients written with respect to the basis $a_1, \ldots, a_{\ell}$. The lemma follows from these observations. \end{proof} 

\tpoint{Main result} We can now state and prove the main result of this paper. 

\begin{theorem} \label{thm:main-fin} 
    For any element\footnote{We could also have worked with elements $x \in \gg^{\tau}_{\pol, \R}$, see \S \ref{subsubNegIwasawa}.} $x \in \wh{B}^- \eta \left( \tau \right)$ with $0 < \tau < 1$, the pro-Hermitian vector bundle $\Psi \left( x \right)$ satisfies the strong summability condition \eqref{strongSummabilityCondition} for any $\varepsilon > 0$, and is hence theta-finite.
\end{theorem}

\begin{proof}
    Consider a scalar $\varepsilon > 0$. We need to prove that the condition \be{}\label{strongSummabilityConditionBundleFromLoop}
        \begin{array}[t]{ccccc}
            \text{\textbf{Sum}} \left( \overline{Y}_{\bullet, x} \otimes \overline{\mathcal{O}} \left( \varepsilon \right) \right) & : & \sideset{}{}\sum\limits_{n \geqslant 0} \; h_{\vartheta}^0 \left( \overline{\ker q_n} \otimes \overline{\mathcal{O}} \left( \varepsilon \right) \right) & < & + \infty
        \end{array}
    \ee is satisfied. Let us write the element $x$ as $x=u^- \, h \, \eta \left( \tau \right)$, with $0 < \tau < 1 $. Recalling the definition $\lambda_1$ for the length of the shortest non-zero element in a lattice (\textit{cf.} \eqref{defn:shortestvector} and the abuse of notation introduced there)
    \be{}\label{shortestLengthInequality} \begin{array}{lllllll} \multicolumn{3}{l}{\lambda_1 \left( \overline{\ker q_n} \otimes \overline{\mathcal{O}} \left( \varepsilon \right), e^{-\varepsilon} \left \Vert \cdot \right \Vert_{V/V^n, x} \right)} & = & \multicolumn{2}{l}{\lambda_1 \left( \overline{\ker q_n} \otimes \overline{\mathcal{O}} \left( \varepsilon \right), e^{-\varepsilon} \left \Vert \cdot \right \Vert_{u^-(0) h \eta \left( \tau \right)} \right)} & \text{ using  \eqref{unip:no-quot}} \\[1em] \qquad \qquad \qquad \qquad & \geqslant & \multicolumn{3}{l}{C_2 e^{-A_2 \sqrt{n}} \lambda_1 \left( \overline{\ker q_n} \otimes \overline{\mathcal{O}} \left( \varepsilon \right), e^{-\varepsilon} \left \Vert \cdot \right \Vert_{h \eta \left( \tau \right)} \right)} && \text{ using  \eqref{unip:lowerbound}} \\[1em] \qquad \qquad \qquad \qquad & \geqslant & \multicolumn{3}{l}{C_1 C_2 e^{- \left( A_1 + A_2 \right) \sqrt{n}} \tau^{-n} \lambda_1 \left( \overline{\ker q_n} \otimes \overline{\mathcal{O}} \left( \varepsilon \right), e^{-\varepsilon} \left \Vert \cdot \right \Vert \right)} && \text{ using  \eqref{hv:lower-bound}} \\[1em] \qquad \qquad \qquad \qquad & \geqslant & \multicolumn{3}{l}{C_1 C_2 e^{- \left( A_1 + A_2 \right) \sqrt{n}} \tau^{-n} e^{- \varepsilon}}, \end{array} \ee the last inequality being a consequence of the integrality of the untwisted norm $\left \Vert \cdot \right \Vert$. We further have \be{}\label{ineqRank} \begin{array}{lllll} \displaystyle \left( \frac{\rk_{\zee} V_{\zee} \left[ n \right]}{2\pi} \right)^{1/2} & = & \displaystyle \left( \frac{\dim_{\C} V \left[ n \right]}{2\pi} \right)^{1/2} & \leqslant & \displaystyle C_4 e^{A_4 \sqrt{n}} \end{array} \ee for two constants $C_4 > 0$ and $A_4 > 0$ independent of $n$. Having chosen $\tau$ such that $0 < \tau < 1$, we have \be{}\label{ineqShortestLengthGr} \begin{array}{lll} \lambda_1 \left( \overline{\ker q_n} \otimes \overline{\mathcal{O}} \left( \varepsilon \right), e^{-\varepsilon} \left \Vert \cdot \right \Vert_{V/V^n, x} \right) & \geqslant & \displaystyle \left( \frac{\rk_{\zee} V_{\zee} \left[ n \right]}{2\pi} \right)^{1/2} \end{array} \ee for $n$ large enough, allowing us to apply \eqref{upperBoundh0} with the estimate \eqref{inequalityConstantCrlambda}. We get \be{}\label{boundH0} \begin{array}{llll} h_{\vartheta}^0 \left( \overline{\ker q_n} \otimes \overline{\mathcal{O}} \left( \varepsilon \right) \right) & \leqslant & \multicolumn{2}{l}{3^{C_4 e^{A_4 \sqrt{n}}} \left( 1 - \frac{C_3}{2 \pi^2 C_1 C_2} e^{\varepsilon} e^{\left( A_1 + A_2 + A_3 \right) \sqrt{n}} \tau^n \right)^{-1}} \\[0.5em] && \qquad \qquad \qquad \qquad & \exp \left( - \pi C_1^2 C_2^2 e^{-2 \varepsilon} e^{-2 \left( A_1 + A_2 \right) \sqrt{n}} \tau^{-2n} \right) \\[1em] & \leqslant & \multicolumn{2}{l}{2 \cdot 3^{C_4 e^{A_4 \sqrt{n}}} \exp \left( - \pi C_1^2 C_2^2 e^{-2 \varepsilon} e^{-2 \left( A_1 + A_2 \right) \sqrt{n}} \tau^{-2n} \right)}  \end{array} \ee for $n$ large enough, as $\tau^n$ decays quicker than $e^{\sqrt{n}}$ grows. Using estimate \eqref{boundH0} and lemma \ref{lem:convergence} below, we note that the strong summability condition \eqref{strongSummabilityConditionBundleFromLoop} is satisfied for any real number $\varepsilon > 0$. We finally obtain theta-finiteness (\textit{cf.} \eqref{thetaFiniteness}) as a consequence.  \end{proof}

\spoint In the proof of Theorem \ref{thm:main-fin}, we used the following elementary convergence lemma whose proof we record here for completeness.

    \begin{lemma}
    \label{lem:convergence}    
        Let $A, B, C, D > 0$ be positive real constants, and assume we have $0 < \tau < 1$. We have
        \begin{equation}
            \begin{array}[t]{lll}
                \sideset{}{}\sum\limits_{n \geqslant 0} \; \exp \left( - A e^{-B \sqrt{n}} \tau^{-2n} + C e^{D \sqrt{n}} \right) & < & + \infty.
            \end{array}
        \end{equation}
    \end{lemma}
    
     \begin{proof}
        We begin by noting that for $0 < \tau < 1$, we have 
        \begin{equation}
        \label{eq:convLemma}
            \begin{array}[t]{lll}
                -A e^{-B \sqrt{n}} \tau^{-n} + C e^{D \sqrt{n}} \tau^n & \underset{n \rightarrow + \infty}{\longrightarrow} & - \infty
            \end{array}
        \end{equation}
        since, in the left-hand side of \eqref{eq:convLemma}, the first term diverges to $- \infty$, while the second term converges to $0$. Thus, there exists a positive integer $n_0$ such that we have
        \begin{equation}
            \begin{array}[t]{llll}
                -A e^{-B \sqrt{n}} \tau^{-n} + C e^{D \sqrt{n}} \tau^n & < & -1 & \text{for } \; n \geqslant n_0.
            \end{array}
        \end{equation}
        Summing over such integers, we get
        \begin{equation}
            \begin{array}[t]{lll}
                \sideset{}{}\sum\limits_{n \geqslant n_0} \; \exp \left( - A e^{-B \sqrt{n}} \tau^{-2n} + C e^{D \sqrt{n}} \right) & \leqslant & \sideset{}{}\sum\limits_{n \geqslant n_0} \; \exp \left( - \tau^{-n} \right).
            \end{array}
        \end{equation}
        Let us now note that the function $x \mapsto \left( 1/x \right) e^{-1/x}$ converges to $0$ as $x \rightarrow 0^+$. Hence, there exists a positive integer $n_1 \geq n_0$ such that we have $\tau^{-n} \exp \left( - \tau^{-n} \right) < 1$ for every $n \geqslant n_1$, thus giving
        \begin{equation}
            \begin{array}[t]{lllll}
                \sideset{}{}\sum\limits_{n \geqslant n_1} \; \exp \left( - \tau^{-n} \right) & \leqslant & \sideset{}{}\sum\limits_{n \geqslant n_1} \; \tau^n & < & + \infty
            \end{array}
        \end{equation}
        using once again the fact that we have $0 < \tau < 1$. This concludes the proof of the lemma.
     \end{proof}

%\nocite{*}
\bibliographystyle{amsplain}
\bibliography{bibliography}

 \end{document}